\documentclass[9pt,shortpaper,twoside,web]{ieeecolor}
\usepackage[table,pdftex,hyperref,svgnames]{xcolor}
\usepackage{generic}
\usepackage{cite}
\usepackage{amsmath,amssymb,amsfonts}
\usepackage{algorithmic}
\usepackage{graphicx}
\usepackage{textcomp}
\usepackage{stmaryrd}
\usepackage{multicol} \usepackage{mathtools}
\usepackage{version,xspace}
\usepackage{url,doi}

\usepackage{hyperref}%
\hypersetup{colorlinks=true, linkcolor=blue, breaklinks=true, urlcolor=blue}
\renewcommand{\theenumi}{(\roman{enumi})}

\newcommand{\ds}{\displaystyle}

\newcounter{saveenum}

 \newcommand{\setdef}[2]{\{#1
	\; | \; #2\}}

\newcommand{\jac}[1]{D\mkern-2.5mu{#1}}
\newcommand{\WP}[2]{\left\llbracket{#1}, {#2}\right\rrbracket}

\newcommand\oprocendsymbol{\hbox{$\triangle$}}
\newcommand\oprocend{\relax\ifmmode\else\unskip\hfill\fi\oprocendsymbol}

                  \let\labelindent\relax
\usepackage{enumitem} 

\newcommand{\verti}[1]{\left\| #1 \right\|}

\DeclareSymbolFont{bbold}{U}{bbold}{m}{n}
\DeclareSymbolFontAlphabet{\mathbbold}{bbold}
\newcommand{\vect}[1]{\mathbbold{#1}}

\newcommand{\vectorzeros}[1][]{\vect{0}_{#1}}

\newcommand{\real}{\mathbb{R}}

\renewcommand{\top}{\mathsf{T}} 

\newtheorem{theorem}{Theorem}[section]
\newtheorem{proposition}[theorem]{Proposition}
\newtheorem{corollary}[theorem]{Corollary}
\newtheorem{lemma}[theorem]{Lemma}
\newtheorem{definition}[theorem]{Definition}
\newtheorem{remark}[theorem]{Remark}
\newtheorem{example}[theorem]{Example}

\newcommand{\suchthat}{\;\ifnum\currentgrouptype=16 \middle\fi|\;}
\newcommand{\scirc}{\raise1pt\hbox{$\,\scriptstyle\circ\,$}}

\def\BibTeX{{\rm B\kern-.05em{\sc i\kern-.025em b}\kern-.08em
    T\kern-.1667em\lower.7ex\hbox{E}\kern-.125emX}}
\begin{document}
\title{Contraction Theory for Positive and Monotone Systems\\ with
  Applications\thanks{This work was supported in part by the Defense Threat
    Reduction Agency under Contract No.~HDTRA1-19-1-0017.}}

\title{Non-Euclidean Contraction Theory for
  Monotone and Positive Systems\thanks{This work was supported in part by the Defense Threat
    Reduction Agency under Contract No.~HDTRA1-19-1-0017.}}

\author{Saber Jafarpour, \IEEEmembership{Member, IEEE}, Alexander Davydov,
  \IEEEmembership{Student Member, IEEE}, and \\ Francesco Bullo,
  \IEEEmembership{Fellow, IEEE} \thanks{Authors are with the Center for
    Control, Dynamical Systems, and Computation, University of California,
    Santa Barbara, 93106-5070, USA. {({\tt \{saber, davydov,
        bullo\}@ucsb.edu})}}}

\maketitle

\begin{abstract}
In this note we study contractivity of monotone
    systems and exponential convergence of positive systems using
    non-Euclidean norms.
  We first introduce the notion of conic matrix measure as a framework
  to study stability of monotone and positive systems.
  We study properties of the conic matrix measures and investigate their connection
  with weak pairings and standard matrix measures.
  Using conic matrix measures and weak pairings, we characterize
  contractivity and incremental stability of monotone systems with
  respect to non-Euclidean norms.
  Moreover, we use conic matrix measures to provide sufficient
  conditions for exponential convergence of positive systems to their
  equilibria.
  We show that our framework leads to novel results on (i) the
  contractivity of excitatory Hopfield neural networks, and (ii) the
  stability of interconnected systems using non-monotone positive
  comparison systems.
\end{abstract}


\section{Introduction}

\paragraph*{Problem description and motivation}

A dynamical system is monotone if its trajectories preserve a partial order
of their initial conditions and is positive if the non-negative orthant is
a forward invariant set. Monotonicity appears naturally in real world
applications including biological systems~\cite{EDS:07}, transportation and
flow networks~\cite{GC-EL-KS:15}, and epidemic networks~\cite{AK-TB-BG:16},
as well as in small-gain analysis of large-scale interconnected
systems~\cite{BSR-CMK-SRW:10,SND-BSR-FRW:10}. Positive systems are also
abundant in engineering and science, for instance, in population
dynamics~\cite{JH-KS:98} and queuing systems~\cite{LF-SR:00}. While the notions of monotonicity and positivity
are identical for linear systems, they are distinct and lead to different
transient and asymptotic behaviors for nonlinear systems. Linear and
nonlinear monotone systems have been studied extensively in dynamical
systems~\cite{HLS:95} and control
theory~\cite{AR:15,EDS:07}. Monotonicity of dynamical systems with respect to arbitrary cones are studied
in~\cite{MWH-HLS:03} and a theory of monotone systems on partially
ordered Banach spaces has been developed in~\cite{AGM:11}.

Contraction theory is a classic
framework~\cite{WL-JJES:98,ZA-EDS:14b,MdB-DF-GR-FS:16,FB:22-CTDS} aimed at
establishing rigorous nonlinear stability properties of dynamical
systems. A dynamical system is contracting if every two
trajectories converge exponentially to one another. Contracting
systems exhibit many desirable asymptotic properties: (i) their asymptotic
behavior is independent of their initial condition, (ii) when the vector
field is time-invariant every trajectory converges to a unique equilibrium
point, and (iii) when the vector field is periodic, every trajectory
converges to a unique periodic orbit. Contracting systems enjoy
also desirable transient behavior and robustness properties including
input-to-state stability in the presence of bounded unmodeled
dynamics. 

While classical approaches mostly focus on contraction with respect to the
$\ell_2$-norm, recent works have shown that stability of monotone and positive system can
be studied more systematically and efficiently using non-Euclidean
norms. It is known that for a monotone system satisfying a conservation law
(resp. translational symmetry), contractivity naturally arises with respect
to $\ell_1$-norms (resp. $\ell_{\infty}$-norms). Contraction of monotone
systems with respect to state-dependent non-Euclidean norms has been
studied in~\cite{SC:19}.  Contraction of monotone systems with respect to
$\ell_1$-norm has been studied for flow networks in~\cite{GC-EL-KS:15}, for
traffic networks in~\cite{GC-EL-KS:15}, and for gene translation systems
in~\cite{MM-EDS-TT:14}. Another relevant topic for monotone systems is the
search for sum-separable and max-separable Lyapunov
functions~\cite{HRF-BB-MJ:18}. Recent works have used contraction with
respect to non-Euclidean norms for monotone systems to find separable
Lyapunov functions~\cite{IRM-JJES:17,YK-BB-MC:20}. Despite all these works,
a differential and integral characterization of monotone and positive
contracting systems with respect to non-Euclidean norms is missing.

\paragraph*{Contribution}
In this note, we build on the framework proposed
in~\cite{AD-SJ-FB:20o} and introduce the notion of conic matrix
measure, characterize its properties, and propose efficient methods
for computing it. We provide a complete characterization of
contractive monotone systems using the one-sided Lipschitz
constant of their vector fields and the conic matrix measure of their
Jacobians. We also propose a sufficient condition, based on the
conic matrix measures, for exponential convergence of positive
systems to equilibrium points.  As a first application of our monotone contraction
framework, we provide
a sufficient condition for contractivity of excitatory Hopfield
neural networks. We remark that strong contractivity of Hopfield
neural networks automatically leads to their global stability for
time-invariant inputs, their entrainment to a unique periodic orbit
for periodic inputs, and their input-to-state stability for general
time-varying inputs. {\color{black} As a second application, we establish a novel
framework for studying input-to-state stability of interconnected
systems.} Our framework is based on comparison with positive dynamical
systems and can accommodate both inhibitory and excitatory
interconnections between subsystems. By allowing the comparison system
to be positive instead of monotone, our framework generalizes the well-known Matrosov-Bellman comparison lemma and
unifies several existing small-gain theorems and comparison lemmas in
the literature.

\section{Notation} 
\paragraph*{Functions, norms and matrix measures}
Let $f:\real_{\ge 0}\to \real_{\ge 0}$ be a function. If $f$ is
differentiable, then we denote its derivative by $f'$. If $f$ is
continuous, we denote its upper {\color{black}Dini derivative} by $D^+f$. We say $f$ is of
class $\mathcal{K}$ if it is strictly increasing and $f(0)=0$. We say $f$
is of class $\mathcal{K}_{\infty}$ if it belongs to class $\mathcal{K}$ and
$\lim_{x\to+\infty} f(x)=\infty$. We say a continuous function
$g:[0,a)\times [0,\infty)\to [0,\infty)$ is of class $\mathcal{KL}$ if,
for each fixed $y$, the map $x\mapsto g(x,y)$ is of class
$\mathcal{K}$ and, for each fixed $x$, the map $y\mapsto g(x,y)$ is
decreasing such that $\lim_{y\to+\infty}g(x,y)=0$. {\color{black} For vectors $v,w\in \real^n$, the Hadamard product of $v$ and
  $w$ is the vector $v\circ w \in \real^n$ define by $(v\circ w)_i = v_iw_i$, for every $i\in
  \{1,\ldots,n\}$.} A matrix $A\in
\real^{n\times n}$ is non-negative if $A_{ij}\ge 0$, for every
$(i,j)\in \{1,\ldots,n\}$, For every matrix $A\in
\real^{n\times n}$, the positive part of $A$ is the matrix
$[A]^{+}\in \real^{n\times n}_{\ge 0}$ defined by $[A]_{ij}^{+}=A_{ij}$ if $A_{ij}\ge 0$ and $[A]^{+}_{ij}=0$
if $A_{ij}<0$. Given $x,y\in \real^n$, 
$x\le y$ if we have $x_i\le y_i$, for every $i\in
\{1,\ldots,n\}$ and we define $[x,y]\subset \real^n$ as the set of all $z\in \real^n$
such that $x\le z\le y$. For a vector $\eta\in \real^n$, the diagonal
matrix $[\eta]\in \real^{n\times n}$ is defined by
$[\eta]_{ii}=\eta_i$, for every $i\in \{1,\ldots,n\}$. Given $A,B\in \real^{n\times n}$, $A\preceq B$ if
$B-A$ is a positive semi-definite matrix. A norm $\|\cdot\|$ on $\real^n$ is
\emph{monotonic}, if for every $x,y\in \real^n$ such that $|x|\le |y|$, we
have $\|x\|\le \|y\|$. {\color{black}For $p\in [1,\infty]$ and $R\in\real^{n\times{n}}$: if
$R$ is positive and diagonal, the $R$-weighted $\ell_p$-norm is a
monotonic norm.} Let $\|\cdot\|$ be a norm on $\real^n$, the
  induced matrix norm on $\real^{n\times n}$ is again denoted by
  $\|\cdot\|$. Given a matrix $A\in
\real^{n\times n}$,  the matrix measure of $A$ with respect to
$\|\cdot\|$ is defined by $\mu(A) := \lim_{h\to
  0^{+}}\frac{\|I_n+hA\|-1}{h}$.

\paragraph*{Weak pairings}
We briefly the notion of a weak pairing (WP) on $\real^{n}$
from~\cite{AD-SJ-FB:20o}. A \emph{WP} on $\real^n$ is a map
$\WP{\cdot}{\cdot}: \real^n \times \real^n \to \real$ satisfying:
\begin{enumerate}
\item\label{WP1}(Sub-additivity and continuity of first argument)
  $\WP{x_1+x_2}{y} \leq \WP{x_1}{y} + \WP{x_2}{y}$, for all
  $x_1,x_2,y \in \real^n$ and $\WP{\cdot}{\cdot}$ is continuous in its
  first argument,
\item\label{WP3}(Weak homogeneity)
  $\WP{\alpha x}{y} = \WP{x}{\alpha y} = \alpha\WP{x}{y}$ and
  $\WP{-x}{-y} = \WP{x}{y}$, for all
  $x,y \in \real^n, \alpha \geq 0$,
\item\label{WP4}(Positive definiteness) $\WP{x}{x} > 0$, for all
  $x \neq \vectorzeros[n],$
\item\label{WP5}(Cauchy-Schwarz inequality) \\
  $|\WP{x}{y}| \leq \WP{x}{x}^{1/2}\WP{y}{y}^{1/2}$, for all
  $x, y \in \real^n.$
\end{enumerate}
For every norm $\|\cdot\|$ on $\real^n$, there exists a (possibly not
unique) associated WP $\WP{\cdot}{\cdot}$ such that $\|x\|^2=\WP{x}{x}$,
for every $x\in \real^n$. A WP $\WP{\cdot}{\cdot}$ satisfies
\emph{Deimling's inequality} if $\WP{x}{y} \le \|y\|\lim_{h\to
  0^{+}}\frac{\|y+hx\|-\|y\|}{h}$, for every $x,y\in \real^n$ and satisfies
the \emph{curve norm derivative formula} if, for every differentiable $x:
(a,b) \to \real^n$ and for almost every $t \in (a,b)$ we have
$\|x(t)\|D^+\|x(t)\| = \WP{\dot{x}(t)}{x(t)}$. For every $p\in (1,\infty)$
and invertible $R\in \real^{n\times n}$, we define
$\WP{\cdot}{\cdot}_{p,R}$ by
\begin{align}\label{eq:p-R}
\WP{x}{y}_{p,R}= \|y\|_{p,R}^{2-p}(Ry \circ |Ry|^{p-2})^\top Rx.
\end{align}
For invertible $R\in \real^{n\times n}$, we define
$\WP{\cdot}{\cdot}_{1,R}$ and $\WP{\cdot}{\cdot}_{\infty,R}$ by
\begin{align}
\WP{x}{y}_{1,R}&= \|Ry\|_1\mathrm{sign}{(Ry)}^{\top}Rx, \label{eq:1-R}\\
\WP{x}{y}_{\infty,R}&= \max_{i\in I_{\infty}(Ry)} (Ry)_i(Rx)_i, \label{eq:inf-R}
\end{align}
where $I_{\infty}(x)=\setdef{i\in \{1,\ldots,n\}}{x_i=\max_i\{|x|_i\}}$. It
can be shown that, for every $p\in [1,\infty]$ and invertible matrix $R\in
\real^{n\times n}$, we have $\|x\|^2_{p,R}=\WP{x}{x}_{p,R}$ and
$\WP{\cdot}{\cdot}_{p,R}$ satisfies Deimling's inequality and the curve
norm derivative formula. We refer to~\cite{AD-SJ-FB:20o} for a detailed
discussion on WPs.

\paragraph*{Dynamical systems}
Consider the dynamical system $\dot{x} = f(t,x)$ on $\real^n$. Let
$\phi(t,t_0,x)$ denote the flow of $f$ at time $t$ starting at time $t_0$
from $x_0$. The vector field $f$ is \emph{positive} if $\real^n_{\ge 0}$ is
a forward invariant set. Let $\mathcal{C}$ be a convex forward invariant
set for vector field $f$. The vector field $f$ is \emph{monotone} on
$\mathcal{C}$, if for every $x_0,y_0\in \mathcal{C}$ such that $x_0\le
y_0$, we have $\phi(t,t_0,x_0)\le \phi(t,t_0,y_0)$, for every $t\ge
t_0$. The Jacobian of $f$ is denoted by $\jac{f}(t,x)$. 
Let $\|\cdot\|$ be a norm with associated WP
$\WP{\cdot}{\cdot}$.The vector field $f$ is \emph{contracting}
with rate $c>0$ if, for $x,y\in \real^n$ and every $t_0\le t\in \real_{\ge 0}$, we
have
\begin{align*}
  \|\phi(t,t_0,x)-\phi(t,t_0,y)\|\le e^{-c(t-t_0)}\|x-y\|.
  \end{align*}
and it is incrementally exponentially stable with rate $c>0$
if, there exists $M>0$ such that, for $x,y\in \real^n$ and every
$t_0\le t \in
\real_{\ge 0}$
\begin{align*}
  \|\phi(t,t_0,x)-\phi(t,t_0,y)\|\le Me^{-c(t-t_0)}\|x-y\|.
  \end{align*}

\section{Conic matrix measures}

  In classical contraction theory, the incremental stability of
  dynamical systems are ensured by imposing suitable conditions on matrix
  measures of their Jacobian.
  Monotonicity or positivity of dynamical systems induces a natural 
  partial order structure on their flows, something which can be used
  to relax the matrix measure conditions for incremental stability.
  In this section, we introduce the notion of conic matrix measures and show how it can be
  used to prove incremental stability of monotone and exponential
  convergence of positive
  systems. 

\begin{definition}[Conic matrix measure]\label{def:cone}
  Let $\|\cdot\|$ be a norm.  The \emph{conic matrix measure} of $A\in
  \real^{n\times n}$, denoted by $\mu^{+}(A)$, is
  \begin{align*}
    \mu^{+}(A) := \lim_{h\to 0^+}\sup_{x\ge \vect{0}_n, x\ne \vect{0}_n}\frac{\|(I_n + h A)x\|/\|x\| -1}{h}.
  \end{align*}
\end{definition} 
Now we study properties of the conic matrix measure. We first
state the following useful lemma.

\begin{lemma}[Monotonicity of $\WP{\cdot}{\cdot}_{p,R}$]\label{lem:useful}
  Let $p\in [1,\infty]$ and $R\in \real^{n\times n}$ be an invertible
  non-negative matrix. Then
  \begin{enumerate}
  \item\label{p1:non-positive} for every $x,y\ge \vect{0}_n$, we have
    $\WP{-x}{y}_{p,R}\le 0$.
  \item\label{p2:inequalityWP} for every $x\le z$ and
    $y\ge \vect{0}_n$, we have $\WP{x}{y}_{p,R}\le \WP{z}{y}_{p,R}$.
  \end{enumerate}
\end{lemma}
\smallskip
\begin{proof}
  Regarding part~\ref{p1:non-positive}, since $R$ is non-negative,
  then we have $Ry\ge \vect{0}_n$ and $Rx\ge \vect{0}_n$. The result
  then follows by using formulas~\eqref{eq:p-R}, \eqref{eq:1-R}, and
  \eqref{eq:inf-R}. Regarding part~\ref{p2:inequalityWP}, note that
  since $x\le z$, there exists $w\ge \vect{0}_n$ such that
  $x = z - w$. Thus, using the subadditivity of the WP, we get
  \begin{align*}
    \WP{x}{y}_{p,R} =  \WP{z-w}{y}_{p,R} \le \WP{z}{y}_{p,R} +
    \WP{-w}{y}_{p,R} \le \WP{z}{y}_{p,R}.
  \end{align*}
  where the last inequality holds because of
  part~\ref{p1:non-positive}.
\end{proof}


  \begin{theorem}[Properties of the conic matrix measure]\label{lem:lemma}
    Let $\|\cdot\|$ be a norm with associated WP $\WP{\cdot}{\cdot}$
    satisfying Deimling's inequality. For every $A,B\in \real^{n\times n}$
    and every $a\in \real$,
    \begin{enumerate}
    \item\label{p5:scalar} $\mu^{+}(aA) = |a|\mu^{+}(\mathrm{sign}(a)A)$;
    \item\label{p6:triangle} $\mu^{+}(A+B)\le \mu^{+}(A)+\mu^{+}(B)$; 
    \item\label{p7:additive} $\mu^{+}(A+a I_n) = \mu^+(A)+a$;
    \item\label{p3:Lumer}
      $\ds \mu^{+}(A) = \sup_{x\ge \vect{0}_n, x\ne \vect{0}_n}
      \frac{\WP{Ax}{x}}{\|x\|^2} $,
    \item\label{p1:inequality}  $\mu^{+}(A)\le \mu(A)$;
        
    \item\label{p2:Metzler} if $\|\cdot\|$ is monotonic and $A$ is
      Metzler then $\mu^{+}(A) = \mu(A)$;
      
     \item\label{p4:FrancescoMajorant} if $\|\cdot\|=\|\cdot\|_{p,R}$
      for $p\in [1,\infty]$ and $R\in \real^{n\times n}$
      invertible and non-negative, then $\mu^+(A) \le \mu^{+}(A+\Delta)$ for every $\Delta\in
      \real^{n\times n}_{\ge 0}$.  
    \end{enumerate}
  \end{theorem}
  \smallskip
  \begin{proof}
    Regarding parts~\ref{p5:scalar}~\ref{p6:triangle},
      and~\ref{p7:additive}, the proofs are straightforward using
      definition of the conic matrix measure in Definition~\ref{def:cone}.

    Regarding~\ref{p3:Lumer}, for every non-zero $x\ge \vect{0}_n$ and
    $h>0$, 
      \begin{align}\label{eq:good-2}
        \|(I_n-h A)x\| &\ge
        \frac{1}{\|x\|}\WP{(I_n-h A)x}{x} \ge (1 -h\frac{\WP{Ax}{x}}{\|x\|^2})\|x\|\nonumber\\ &\ge \left(1-h\sup_{x\ge
        \vect{0}_n,x\ne \vect{0}_n}\frac{\WP{Ax}{x}}{\|x\|^2}\right)\|x\|,
      \end{align}
      where first inequality holds by Cauchy-Schwarz, the second inequality is by
      subadditivity of WP and the fact that $h>0$, and the
      third inequality holds because $h>0$.
      Now consider $F(h) = (I_n-h A)^{-1}$. By some simple
      algebraic manipulation, we get
      \begin{align}\label{eq:algebraic}
        F(h) = I_n + h  A + h^2 A^2 F(h).
      \end{align}
      Note that for
      $h=0$, we have $F(0)=I_n$ and thus, for every non-zero $x\ge 0$, we
      have $F(h)x = (I_n-h A)^{-1}x = x \ge \vect{0}_n$. Since
      $\mathbb{S}_{\|\cdot\|}^n=\setdef{x\in \real^n}{\|x\|=1}$ is a compact
      set, using a continuity argument, there exists small enough $h^*>0$ such that $(I_n-h
      A)^{-1}x\ge \vect{0}_n$ for every $x\in \mathbb{S}_{\|\cdot\|}^n$ and $x\ge
      \vect{0}_n$ and every $0 \le h\le h^*$. Thus, for every
      $0\le h\le h^*$ and every $v\ge \vect{0}_n$ such that $v\ne \vect{0}_n$, we have
      \begin{align*}
        \|F(h)v\|/\|v\| &=\frac{\|x\|}{\|(I_n-hA)x\|} \le \frac{1}{(1-h\sup_{\|x\|=1,x\ge
        \vect{0}_n}\frac{\WP{Ax}{x}}{\|x\|^2})}, \\
        h^2\|A^2(&I_n-h A)^{-1}v\|/\|v\| \le  h^2\|A^2(I_n-h A)^{-1}\|,
      \end{align*}
      where the first equality holds by the change of coordinate
      $x = (I_n-h A)^{-1}v$ and
      the second inequality holds by~\eqref{eq:good-2}.
      Therefore, for every $0\le h\le h^*$, we have
      \begin{align*}
        \mu^{+}&(A) = \lim_{h\to 0^{+}}\sup_{\|v\|=1,v\ge \vect{0}_n}\frac{\|(I_n+h
                 A)v\|/\|v\|-1}{h} \\ &\le \lim_{h\to 0^{+}}\sup_{\|v\|=1,v\ge \vect{0}_n}\frac{\|F(h)v\|+h^2\|A^2F(h)v\|-\|v\|}{h\|v\|}
        \\ & = \lim_{h\to 0^{+}}\sup_{\|v\|=1,v\ge
             \vect{0}_n}\frac{\|F(h)v\|-\|v\|}{h\|v\|} \\& \le \lim_{h\to 0^{+}}\frac{1}{h} \Big(\frac{1}{1
                                                            -h\sup_{\|x\|=1,x\ge
                                                            \vect{0}_n}\frac{\WP{Ax}{x}}{\|x\|^2}}-1\Big)
                                                             \\ &=
                                                            \sup_{\|x\|=1,x\ge
                                                            \vect{0}_n}\frac{\WP{Ax}{x}}{\|x\|^2}
                                                                  = \sup_{x\ge
                                                            \vect{0}_n,
                                                                  x\ne
                                                                  \vect{0}_n}\frac{\WP{Ax}{x}}{\|x\|^2},
      \end{align*}
      where the first equality holds by definition, the second inequality holds by applying triangle
      inequality to the algebraic
                 equation~\eqref{eq:algebraic}, and the fourth
                 inequality holds by~\eqref{eq:good-2}.  
      This means that $\mu^+(A)\le \sup_{x\ne
        \vect{0}_n, x\ge \vect{0}_n}\frac{\WP{Ax}{x}}{\|x\|^2}$. Additionally, using Deimling's inequality,
      \begin{align*}
        &\WP{Ax}{x}\le \|x\|\lim_{h\to
        0^{+}}\frac{\|x+hAx\|-\|x\|}{h} \\ & \le  \|x\|^2 \lim_{h\to
        0^{+}}\sup_{x\ge \vect{0}_n, x\ne \vect{0}_n}\frac{\|x+hAx\|/\|x\|-1}{h}=\|x\|^2\mu^{+}(A).
      \end{align*}
      This means that $\sup_{x\ge \vect{0}_n, x\ne \vect{0}_n}\frac{\WP{Ax}{x}}{\|x\|^2}\le
      \mu^+(A)$ and completes the proof
      of~\ref{p3:Lumer}. Regarding~\ref{p1:inequality}, the proof is
      straightforward using the definitions. Regarding
      part~\ref{p2:Metzler}, we show that, if $\|\cdot\|$ is monotonic,
      then, for small enough $h>0$, $\sup_{x\ge \vect{0}_n}\|(I_n + h
      A)x\|/\|x\| = \sup_{x\ne \vect{0}_n}\|(I_n + h A)x\|/\|x\|$. First
      note that, by definition, $\sup_{x\ge \vect{0}_n}\|(I_n + h
      A)x\|/\|x\| \le \sup_{x\ne \vect{0}_n}\|(I_n + h A)x\|/\|x\|$. Our
      goal in this part is to show the other side of the inequality. Since
      $A$ is Metzler, for small enough $h>0$, the matrix $I_n + hA$ is
      non-negative. By triangle inequality, this implies that $|(I_n+hA)x|
      \le (I_n+hA)|x|$, for every $x\in \real^n$ and for small enough
      $h>0$. Since $\|\cdot\|$ is a monotonic norm, for every $x\in
      \real^n$ such that $x\ne \vect{0}_n$,
    \begin{align*}
      \|(I_n + h A)x\|/\|x\|  \le  \|(I_n + h A)|x|\|/\||x|\|. 
    \end{align*}
    and thus, for every $x\ne \vect{0}_n$ and every small
    enough $h>0$, we have
    \begin{align*}
      \tfrac{\|(I_n + h A)x\|/\|x\|-1}{h}  \le  \tfrac{\|(I_n + h A)|x|\|/\||x|\|-1}{h}. 
    \end{align*}
    Since $|x|\ge \vect{0}_n$, we can define $y=|x|$ and take the sup of both
    sides of the above inequality over $x\ne \vect{0}_n$, 
    \begin{align*}
      \mu(A) & = \lim_{h\to 0^+} \sup_{x\ne \vect{0}_n}\frac{\|(I_n + h
        A)x\|/\|x\|-1}{h}  \\ & \le  \lim_{h\to
        0^+}\sup_{y\ge\vect{0}_n, y\ne \vect{0}_n}\frac{\|(I_n + h A)y\|/\|y\|-1}{h} = \mu^{+}(A).
    \end{align*}
    Regarding part~\ref{p4:FrancescoMajorant}, by
    Lemma~\ref{lem:useful}\ref{p2:inequalityWP}, we have
    \begin{align*}
      \mu^{+}_{p,R}(A) &= \sup_{x\ge
        \vect{0}_n, x\ne
                         \vect{0}_n}\frac{\WP{Ax}{x}_{p,R}}{\|x\|_{p,R}^2}
      \\ & \le \sup_{x\ge
        \vect{0}_n, x\ne \vect{0}_n}\frac{\WP{(A+\Delta)x}{x}_{p,R}}{\|x\|_{p,R}^2}
      = \mu_{p,R}^{+}(A+\Delta).
    \end{align*}\end{proof}
  
  Next, we provide formulas for some useful conic matrix measures.
  
  \begin{theorem}[Computing conic matrix measure]\label{thm:computational}
    Let $A\in \real^{n\times n}$ be a matrix, $\|\cdot\|$ be a norm with the conic
    matrix measure $\mu^+$, $R\in \real^{n\times n}$ be an invertible
    non-negative matrix, and $\eta\in \real^{n}_{>0}$. Then
    \begin{enumerate}
      \item\label{p0:non-negative-measureweights} $\mu^{+}_{R}(A) \le \mu^{+}(RAR^{-1})$;
    \item\label{p1:diagonal} $\mu^{+}_{[\eta]}(A) =
      \mu^{+}([\eta]A[\eta]^{-1})$;

     \item\label{p1.5:conic-2}{\color{black} $\mu^{+}_2(A)=p^*- |\mu_2(A)|$, where $p^*$ is the
         optimal value of following Quadratically
         Constrained Quadratic Program (QCQP):
         \begin{align}\label{eq:qcqp}
         p^* = \max &\;\;x^{\top}(\tfrac{1}{2}(A+A^{\top})+|\mu_2(A)| I_n )x \nonumber\\
                    &\;\; x^{\top} x \le 1, \;\; x\ge \vect{0}_n.
         \end{align}} 
    \item\label{p2:conic-1} $\mu^{+}_1(A) = \max_{j} \{a_{jj} +
      \sum\nolimits_{i\ne j}[a_{ij}]^{+}\}$;
    \item\label{p3:conic-inf} $\mu^{+}_{\infty}(A) = \max_{i}
      \{a_{ii} + \sum\nolimits_{j\ne i}[a_{ij}]^{+}\}$.
    \end{enumerate}
    Moreover, if $A$ is Metzler, then the following statements hold:
    \begin{enumerate}\setcounter{enumi}{5}
    \item\label{p1:l1-identity}
      $\mu^{+}_{1,[\eta]}(A)=\mu_{1,[\eta]}(A) =\min\setdef{c\in \real}{\eta^{\top}A \le c\eta^{\top}} $;
    \item\label{p2:linf-identity} $\mu^{+}_{\infty,[\eta]^{-1}}(A) = 
      \mu_{\infty,[\eta]^{-1}}(A) = 
      \min\setdef{c\in \real}{A\eta \le c\eta}$; 
    \item\label{p3:l2-identity} $\mu^{+}_{2,[\eta]}(A) = \mu_{2,[\eta]}(A) =
      \min\setdef{c\in \real}{[\eta]A +
        A^{\top}[\eta]\preceq 2c [\eta]}$.
    \end{enumerate}
  \end{theorem}
  \smallskip
  \begin{proof}
    Regarding part~\ref{p0:non-negative-measureweights},
    we compute
    \begin{align}\label{eq:weighted-conicmeasure}
      &\mu^{+}_{R}(A)= \lim_{h\to 0^{+}}\sup_{x\ge
        \vect{0}_n, x\ne \vect{0}_n}\tfrac{\|(R+hRA)x\|/\|Rx\|-1}{h} \nonumber\\ & \le \lim_{h\to 0^{+}}\sup_{y\ge
        \vect{0}_n}\tfrac{\|(I_n+hRAR^{-1})y\|/\|y\|-1}{h}= \mu^{+}(RAR^{-1}),
    \end{align}
    where the second inequality holds by setting $y = Rx$ and noting
    that if $x\ge \vect{0}_n$, then $y = Rx\ge \vect{0}_n$.  Regarding
    part~\ref{p1:diagonal}, if $R=[\eta]$, then in
    equation~\eqref{eq:weighted-conicmeasure}, we have
    $\sup_{x\ge \vect{0}_n, x\ne \vect{0}_n}\|(R+hRA)x\|/\|Rx\| =
    \sup_{y\ge \vect{0}_n, y\ne \vect{0}_n}\|(I_n+hRAR^{-1})y\|/\|y\|$, since
    $x\ge \vect{0}_n$ if and only if $y\ge \vect{0}_n$. Thus, all the relations in
    equation~\eqref{eq:weighted-conicmeasure} are equality and the
    result follows.  {\color{black} Regarding part~\ref{p1.5:conic-2}, for a symmetric
      matrix $M\in \real^{n\times n}$, we denote the the largest
      eigenvalue of $M$ by $\lambda_{\max}(M)$ and we define $B =
      \tfrac{1}{2}(A+A^{\top})+|\mu_2(A)| I_n$. Therefore, 
      \begin{align*}
        \lambda_{\max}(B) &=
        \lambda_{\max}(\tfrac{1}{2}(A+A^{\top})+|\mu_2(A)| I_n) \\ & = \mu_2(A)
        + |\mu_2(A)| \ge 0 
      \end{align*}
      where the first equality holds by the formula $\mu_2(A) =
      \lambda_{\max}(\frac{1}{2}(A+A^{\top}))$ and last inequality
      holds because $\mu_2(A)\le |\mu_2(A)|$. As
      a result, the matrix $B$ is positive semi-definite. Using
  Theorem~\ref{lem:lemma}\ref{p3:Lumer},
\begin{align*}
      \mu^{+}_{2}(A) = \sup_{x\ge \vect{0}_n, x\ne
  \vect{0}_n}\tfrac{x^{\top}Ax}{\|x\|^2_{2}} =\sup_{y\ge
  \vect{0}_n,\|y\|_2=1}y^{\top}Ay, 
\end{align*}
where the last equality holds by the change of variable
$y=\frac{x}{\|x\|_2}$. This implies that,
\begin{align*}
  \mu^{+}_2(A) &= \mu^{+}_2(\tfrac{1}{2}(A+A^{\top}))
= \mu^{+}_2(B-|\mu_2(A)| I_n) \\ &= \mu_2^{+}(B) - |\mu_2(A)|, 
\end{align*}
where the last equality holds by Theorem~\ref{lem:lemma}\ref{p7:additive}. Moreover, 
\begin{align*}
      \mu^{+}_{2}(B) =\max_{y\ge \vect{0}_n,\|y\|_2=1}y^{\top}By =
  \max_{y\ge \vect{0}_n,\|y\|_2\le 1}y^{\top}By
\end{align*}
where the last equality holds because $B$ is positive semi-definite and
thus $x\mapsto x^{\top}Bx$ is convex. Regarding part~\ref{p2:conic-1}, using
    Theorem~\ref{lem:lemma}\ref{p3:Lumer}, 
    \begin{align*}
      \mu^{+}_{1}(A) &= \sup_{x\ge
        \vect{0}_n, x\ne
      \vect{0}_n}\tfrac{\mathrm{sign}(x)^{\top}Ax}{\|x\|_{1}}
      =\sup_{\|y\|_1=1, y\ge \vect{0}_n}\mathrm{sign}(y)^{\top}Ay \\ &
                                                                       =
      \sup_{\|y\|_1=1, y\ge \vect{0}_n}y^{\top}A^{\top}\mathrm{sign}(y), 
    \end{align*}
    where the second equality holds by the change of variable
    $y=\frac{x}{\|x\|}$ and the last equality holds because
    $\mathrm{sign}(y)^{\top}Ay = y^{\top}A^{\top}
    \mathrm{sign}(y)$. Note that, for every $i\in \{1,\ldots,n\}$ such
    that $y_i\ne 0$, we have $[A^{\top}\mathrm{sign}(y)]_i \le a_{ii}
    + \sum_{k\ne i}[a_{ki}]^+$. This implies that
    \begin{align*}
      \sup_{\|y\|_1=1, y\ge \vect{0}_n}y^{\top}A^{\top}\mathrm{sign}(y) &\le 
      \sup_{\|y\|_1=1, y\ge \vect{0}_n} \sum_{i=1}^{n} y_i (a_{ii} +
      \sum\nolimits_{i\ne k}[a_{ki}]^{+}) \\ & \le \max_{j}\{a_{jj} +
      \sum\nolimits_{i\ne j}[a_{ij}]^{+}\}.
    \end{align*}
    Now we show the converse inequality, Assume that $k\in \{1,\ldots,n\}$ is such that
    \begin{align*}
     a_{kk} +
      \sum\nolimits_{i\ne k}[a_{ik}]^{+} = \max_{j}\{a_{jj} +
      \sum\nolimits_{i\ne j}[a_{ij}]^{+}\}
    \end{align*}
    For every $\epsilon\in (0,1)$, we define $z(\epsilon)\in \real^{n}_{\ge 0}$ as
    follows:
    \begin{align*}
      [z(\epsilon)]_i=\begin{cases}
        1-\epsilon &  i = k,\\
        \epsilon & a_{ik}>0,\\
        0 & a_{ik} = 0
      \end{cases}
    \end{align*}
    Then by setting $y(\epsilon) =
    \frac{z(\epsilon)}{\|z(\epsilon)\|_1}$, we get 
 \begin{multline*}
      \sup_{\|y\|_1=1, y\ge
   \vect{0}_n}yA^{\top}\mathrm{sign}(y)\ge \lim_{\epsilon\to 0^{+}} y(\epsilon)^{\top}A^{\top}\mathrm{sign}(y(\epsilon)) \\=
   \lim_{\epsilon\to 0^{+}} y(\epsilon)_k(a_{kk} +
      \sum\nolimits_{i\ne k}[a_{ik}]^{+}) + \lim_{\epsilon\to 0^{+}}\sum_{i\ne
   k} y(\epsilon)_i[A^{\top}\mathrm{sign}(y)]_i \\ = a_{kk} +
      \sum\nolimits_{i\ne k}[a_{ik}]^{+},
    \end{multline*}
    where the last equality holds because $\lim_{\epsilon\to
      0^{+}}y(\epsilon)_i=1$ and $\lim_{\epsilon\to
      0^{+}}y(\epsilon)_k=0$, for every $k\ne i$. This conclude the proof
    of part~\ref{p2:conic-1}. Regarding
    part~\ref{p3:conic-inf}, using
    Theorem~\ref{lem:lemma}\ref{p3:Lumer}, 
    \begin{align*}
      \mu^{+}_{\infty}(A) & = \sup_{x\ge
        \vect{0}_n, x\ne \vect{0}_n}\max_{i\in
      I_{\infty}(x)}\tfrac{(Ax)_i}{x_i} = \sup_{y\ge
        \vect{0}_n, \|y\|_{\infty}=1}\max_{i\in
      I_{\infty}(y)}(Ay)_i \\ & = \sup_{y\ge
        \vect{0}_n, \|y\|_{\infty}=1}\max_{i\in
      I_{\infty}(y)} \sum_{j=1}^{n} a_{ij}y_j \\ & \leq \max_{i} \{a_{ii} +
      \sum\nolimits_{j\ne i}[a_{ij}]^{+}\}.                  
    \end{align*}
    Now we show the converse inequality. Let $k\in \{1,\ldots,n\}$ be such that
\begin{align*}
     a_{kk} +
      \sum\nolimits_{i\ne k}[a_{ki}]^{+} = \max_{i}\{a_{ii} +
      \sum\nolimits_{j\ne i}[a_{ij}]^{+}\}.
    \end{align*}
    We define $z\in \real^{n}_{\ge 0}$ as follows:
    \begin{align*}
      z_i =\begin{cases}
        1 & i=k \mbox{ or } a_{ki}>0,\\
        0 & a_{ki}\le 0,
        \end{cases}
    \end{align*}
    Then it is easy to see that
    \begin{align*}
      \sup_{y\ge
        \vect{0}_n, \|y\|_{\infty}=1}\max_{i\in
      I_{\infty}(y)}(Ay)_i &\ge \max_{i\in
      I_{\infty}(z)}(Az)_i \\ & \ge (Az)_k = a_{kk} +
      \sum\nolimits_{i\ne k}[a_{ki}]^{+} \\ & = \max_{i}\{a_{ii} +
      \sum\nolimits_{j\ne i}[a_{ij}]^{+}\}.
      \end{align*}
This completes the proof of part~\ref{p3:conic-inf}.}
    Regarding~\ref{p1:l1-identity}, the fact that
    $\mu^{+}_{1,[\eta]}(A)=\mu_{1,[\eta]}(A)$ follows from
    Theorem~\ref{lem:lemma}\ref{p2:Metzler}. Let $b\in \real$ be such that
    $\mu^{+}_{1,[\eta]}(A)\le c$. By part~\ref{p2:conic-1},
       \begin{align*}
         \max(\eta^{\top}A[\eta]^{-1}) = \max(\vect{1}_n^{\top}[\eta]A[\eta]^{-1})
         =\mu^{+}_{1,[\eta]}(A) \le c,
         \end{align*}
       which is also equivalent to $\eta^{\top}A\le
       \eta^{\top}c$. Therefore $\mu^{+}_{1,[\eta]}(A)\le c$ is
       equivalent to $\eta^{\top}A\le \eta^{\top}c$ and the proof is
       complete by taking the $\min$ over $c$. The proof of
       parts~\ref{p2:linf-identity} and~\ref{p3:l2-identity} are similar.      
     \end{proof}

     \begin{remark}{\color{black} The following remarks are in order.}
       \begin{enumerate}
       \item {\color{black} (Computing conic matrix
           measures) Theorem~\ref{thm:computational}\ref{p1.5:conic-2} presents a
           QCQP optimization problem for computing the $\ell_2$-norm conic
           matrix measure. Note that, while the QCQP~\eqref{eq:qcqp} is not in
           general convex, one can find the following convex
           relaxation of~\eqref{eq:qcqp} based on semidefinite
           programming~\cite[Equation (16)]{LV-SB:96}:
           \begin{align}\label{eq:semidef}
             \min_{X\in \real^{n\times n},x\in \real^n} &\quad
                                                          \mathrm{tr}(-(\tfrac{1}{2}(A+A^{\top})+|\mu_2(A)|
                                                          I_n)X)\nonumber\\
                  & \quad \mathrm{tr}(X)\le 1,\;\;\; -x\le \vect{0}_n, \;\;\; \begin{bmatrix}X & x\\ x^\top &
                    1 \end{bmatrix}\succeq 0, 
           \end{align}
           where $\mathrm{tr}$ is the trace operator. The semidefinite
           programming~\eqref{eq:semidef} can be solved using
           CVX~\cite{MG-SB:11-cvx}. Moreover, the semidefinite
           program~\eqref{eq:semidef} provides an upper bound for the optimal solution
           of~\eqref{eq:qcqp}. }
         \item  (Spectral abscissa and conic matrix
         measures) As is shown in Theorem~\ref{lem:lemma}, the conic
         matrix measure shares several nice features with the matrix
         measure including, positive homogeneity, subadditivity, and
         translation properties. Remarkably, unlike the matrix
         measure, the conic matrix measure is sometimes smaller than
         the spectral abscissa. For instance, consider the matrix
         $A=\left[\begin{smallmatrix}-1 &-\frac{1}{2}\\ -1 &
             -2 \end{smallmatrix}\right]$ with eigenvalues
         $\lambda_1(A)=-0.6340$ and $\lambda_2(A)=-2.3660$. Using the
         formula in Theorem~\ref{thm:computational}\ref{p3:conic-inf},
         \begin{align*}
           \mu_{\infty}^{+}(A) = -1 < \max_i
           \mathrm{Re}(\lambda_i(A))< \mu_{\infty}(A) = -0.5. 
           \end{align*} 
         \end{enumerate}
         \end{remark}
       \smallskip
       \smallskip

  
  Finally, we establish a generalized version of Coppel's inequality.

  \begin{theorem}[Conic Coppel's inequality]\label{thm:cone-coppel}
    Let $\|\cdot\|$ be a norm and $t\mapsto A(t)$ be a continuous
    map. Consider the dynamical system
    \begin{align}\label{eq:coppel}
      \dot{x} = A(t) x 
    \end{align}
    If $A(t)$ is Metzler for all $t\ge 0$ and $x(0)\geq\vect{0}_n$, then
    \begin{align*}
      \|x(t)\| \le
      \exp\Big(\int_{0}^{t}\mu^{+}(A(\tau))d\tau\Big)\|x(0)\|,\qquad\mbox{for all
      }t \ge 0.
    \end{align*}
  \end{theorem}
  \smallskip
  \begin{proof}
    Note that
    \begin{align*}
      x(t+h) = x(t) + h A(t) x(t) + \mathcal{O}(h^2) = (I_n + h A(t) ) x(t) + \mathcal{O}(h^2)
    \end{align*}
    and, in turn,
    $D^{+}\|x(t)\| = \lim_{h\to 0^{+}}\frac{\left\|(I_n +
        hA(t))x(t)\right\| - \|x(t)\|}{h}$. Since $A(t)$ is Metzler
    for every $t\ge 0$ and $x(0)\geq\vect{0}_n$, it is well known that
    $x(t) \ge \vect{0}_n$, for every $t\ge 0$. Therefore, for every
    $t\in \real_{\ge 0}$, we obtain
    \begin{align*}
     \lim_{h\to
      0^{+}}\frac{\left\|(I_n + hA(t))x(t)\right\| - \|x(t)\|}{h} \le 
      \mu^{+}(A(t))\|x(t)\|.
    \end{align*}
    The result then follows from Gr\"{o}nwall\textendash{}Bellman Lemma.
  \end{proof}

\section{Contracting monotone and positive systems}

In this section, we use the notions of conic matrix
  measure and WP to study contractive
  monotone systems and converging positive systems. Our first result presents a characterization of contracting
monotone systems using conic matrix measures and WPs.
     
     \begin{theorem}[Contracting monotone systems]\label{thm:contraction-monotone-general}
       Let $\dot{x}=f(t,x)$ be a monotone dynamical system with a
       convex forward invariant set $\mathcal{C}\subseteq\real^n$ and
       $\|\cdot\|$ be a norm with associated WP $\WP{\cdot}{\cdot}$
       satisfying Deimling's inequality. If $\|\cdot\|$ is monotonic,
       then the following statements are equivalent for $b\in \real$:
       \begin{enumerate}
         \item\label{p1:mu-general-1} $\mu^{+}(\jac{f}(t,x))\le b$, for every
        $(t,x)\in \real_{\ge 0}\times \mathcal{C}$;
      \item\label{p3:integral-general-1} $\WP{f(t,x)-f(t,y)}{x-y}\le b\|x-y\|^2$,
        for every $(t,x),(t,y)\in \real_{\ge 0}\times\mathcal{C}$ such
        that $x\ge y$;
        \item \label{p5:contraction-1}
        $\|\phi(t,t_0,x_0)-\phi(t,t_0,y_0)\| \le
        e^{b(t-s)}\|\phi(s,t_0,x_0)-\phi(s,t_0,y_0)\|$, for every $t_0\le s\le t$ and every 
        $x_0,y_0\in \mathcal{C}$.
         \setcounter{saveenum}{\value{enumi}}       
      \end{enumerate}
      Instead, if $\|\cdot\|$ is not monotonic, then
      conditions~\ref{p1:mu-general-1}, 
      \ref{p3:integral-general-1} are equivalent
      and, with $\mathcal{C}=[x_{\min},x_{\max}]$ for
      some $x_{\min}<x_{\max}$, they imply
      \begin{enumerate} \setcounter{enumi}{\value{saveenum}}
      \item \label{p6:contraction+nonmonotone} there exists $M>0$ such
        that, for  every $x_0,y_0\in [x_{\min},x_{\max}]$ and
        every $t_0\le s\le t$,
        \begin{multline*}
      \|\phi(t,t_0,x_0)-\phi(t,t_0,y_0)\| \\ \le M
          e^{b(t-s)}\|\phi(s,t_0,x_0)-\phi(s,t_0,y_0)\|.
        \end{multline*}
        \end{enumerate}
      \end{theorem}
      \smallskip
     \begin{proof}
       Regarding~\ref{p1:mu-general-1}$\implies$\ref{p3:integral-general-1}, compute
      \begin{multline*}
        \WP{f(t,x)-f(t,y)}{x-y}\\ = \WP{\Big(\int_{0}^{1} Df(t,\tau x +
        (1-\tau)y)d\tau\Big)(x-y)}{x-y} \\ \le \int_{0}^{1} \WP{Df(t,\tau x +
        (1-\tau)y)(x-y)}{x-y} d\tau \\ \le \mu^+(Df(t,\tau x +
        (1-\tau)y))\|x-y\|^2\le b\|x-y\|^2,
      \end{multline*}
      where the first equality is by the Mean Value Theorem, the second
      inequality is by the subadditivity of the WP and the third inequality holds by
      Theorem~\ref{lem:lemma}\ref{p3:Lumer} and the fact that
      $x-y\ge\vect{0}_n$.
      Regarding~\ref{p3:integral-general-1}$\implies$\ref{p1:mu-general-1},
      pick $x = y + hv$, for $v\in \real^n_{\ge 0}$ and $h>0$. Thus, 
      \begin{align*}
        &\WP{f(t,x)-f(t,y)}{x-y} = \WP{f(t,y+hv)-f(t,y)}{hv} \\ & = h^2\WP{\tfrac{f(t,y+hv)-f(t,y)}{h}}{v} \le
        b\|x-y\|^2 = bh^2\|v\|^2.
      \end{align*}
      In the limit as $h\to 0^{+}$, for every $y\in \real^n$
      and every $v\in \real^n_{\ge 0}$,
      \begin{align*}
        \WP{\jac{f}(t,y)v}{v} = \lim_{h\to 0^+}\WP{\tfrac{f(t,y+hv)-f(t,y)}{h}}{v} \le
        b\|v\|^2. 
      \end{align*}
      where the first equality holds by the continuity of WP in the first
      argument. As a result, by Theorem~\ref{lem:lemma}\ref{p3:Lumer}, 
      $\mu^{+}(\jac{f}(t,x))\le b$, for every $t\in \real_{\ge 0}$ and
      every $x\in
      \mathcal{C}$. Regarding~\ref{p1:mu-general-1}$\implies$\ref{p5:contraction-1},
      since $f$ is monotone, $Df(t,x)$ is Metzler for
      every $(t,x)\in \real_{\ge 0}\times \mathcal{C}$. Since
      $\|\cdot\|$ is monotonic,
      Theorem~\ref{lem:lemma}\ref{p2:Metzler} implies 
      $\mu^{+}(Df(t,x))=\mu(Df(t,x))$, for every $x\in \mathcal{C}$
      and every $t\ge
      0$. These conclusions then follow from~\cite[Theorem 29]{AD-SJ-FB:20o}.
      Regarding~\ref{p5:contraction-1}$\implies$\ref{p3:integral-general-1},
      note that
      $\|\phi(t+h,t_0,x_0)-\phi(t+h,t_0,y_0)\|\le
      e^{bh}\|\phi(t,t_0,x_0)-\phi(t,t_0,y_0)\|$, for every $h>0$. As a
      result
      \begin{multline*}
        \lim_{h\to
        0^{+}}\tfrac{\|\phi(t+h,t_0,x_0)-\phi(t+h,t_0,y_0)\|-\|\phi(t,t_0,x_0)-\phi(t,t_0,y_0)\|}{h}
        \\ \le \lim_{h\to
             0^+}\frac{e^{bh}-1}{h}\|\phi(t,t_0,x_0)-\phi(t,t_0,y_0)\| \\= b \|\phi(t,t_0,x_0)-\phi(t,t_0,y_0)\|.
      \end{multline*}
      Thus, by Deimling's inequality, for every
      $x_0,y_0\in \mathcal{C}$,
      \begin{multline*}
        \WP{f(t,\phi(t,t_0,x_0))-f(t,\phi(t,t_0,y_0))}{\phi(t,t_0,x_0)-\phi(t,t_0,y_0)}
        \\ \le b \|\phi(t,t_0,x_0)-\phi(t,t_0,y_0)\|^2. 
      \end{multline*}
      This concludes the proof of
      \ref{p5:contraction-1}$\implies$\ref{p3:integral-general-1}. Regarding~\ref{p6:contraction+nonmonotone},
      first we show that, for $x_0\ge y_0$, 
      $\|\phi(t,t_0,x_0)-\phi(t,t_0,y_0)\|\le e^{b(t-s)}\|\phi(s,t_0,x_0)-\phi(s,t_0,y_0)\|$. For
      $\alpha\in[0,1]$, define
      $\psi(t,\alpha)=\phi\big(t,s,\alpha \phi(s,t_0,x_0)+(1-\alpha)\phi(s,t_0,y_0)\big)$ and
      note $\psi(t_0,\alpha) = \alpha \phi(s,t_0,x_0)+(1-\alpha)\phi(s,t_0,y_0)$ and
      $\frac{\partial \psi}{\partial \alpha}(s,\alpha)= \phi(s,t_0,x_0)-\phi(s,t_0,y_0)$.  We
      then compute:
  \begin{align*}
    \frac{\partial}{\partial t} 
    \frac{\partial}{\partial \alpha} \psi(t,\alpha) &= 
    \frac{\partial}{\partial \alpha}   \frac{\partial}{\partial t}
    \psi(t,\alpha) = 
    \frac{\partial}{\partial \alpha}   f(t, \psi(t,\alpha) )
    \\ & = \frac{\partial f}{\partial x} (t, \psi(t,\alpha) ) 
    \frac{\partial}{\partial \alpha}   \psi(t,\alpha) .
  \end{align*}
  Therefore, $\frac{\partial \psi}{\partial \alpha}(t,\alpha)$ satisfies
  the linear time-varying differential equation $ \frac{\partial}{\partial
    t} \frac{\partial \psi }{\partial \alpha} = \jac{f}(t,\psi)
  \frac{\partial \psi}{\partial \alpha} $. Moreover, $x_0-y_0\ge 0$
  and $Df(t,x)$ is Metzler, for every $(t,x)\in \real_{\ge 0}\times
  \mathcal{C}$. Therefore, Theorem~\ref{thm:cone-coppel} implies
  \begin{align}
    \label{eq:bound-Coppel}
    \verti{\tfrac{\partial \psi}{\partial \alpha}(t,\alpha)} &\leq 
   \verti{\tfrac{\partial \psi}{\partial \alpha}(s,\alpha)}\exp\Big(
                                                               \int_{s}^t\mu^{+}\big(
                                                               \jac{f}(t,\psi(\tau,\alpha))
                                                               \big)
                                                               d\tau\Big)\nonumber
                                                               \\ & \leq e^{b(t-s)} \verti{\phi(s,t_0,x_0)-\phi(s,t_0,y_0)},
  \end{align}
  where we used $\mu^{+}(\jac f(t,x))\le b$, for every $t\in
  \real_{\ge 0}$ and  $x\in \real^n$. In turn, 
  inequality~\eqref{eq:bound-Coppel} implies
  \begin{align*}
    \verti{\phi(t,t_0,x_0)-\phi(t,t_0,y_0)}& =\verti{\psi(t,1)-\psi(t,0)}
    \\ &=
    \verti{\int_{0}^{1}\tfrac{\partial \psi(t,\alpha)}{\partial
        \alpha} d\alpha}  \le \int_{0}^{1}\verti{\tfrac{\partial \psi(t,\alpha)}{\partial
        \alpha}}
    d  \alpha \\ &\leq e^{b(t-s)} \verti{\phi(s,x_0)-\phi(s,y_0)}.                            
  \end{align*}
  Now assume that $x_0\not\ge y_0$. We define $\xi,\eta\in \real^n$ by
  \begin{align*}
    \xi_i = \max\{(x_0)_i,(y_0)_i\}, \quad \eta_i = \min\{(x_0)_i,(y_0)_i\}, \quad i\in \{1,\ldots,n\}
  \end{align*}
  Then $\xi,\eta\in [x_{\min},x_{\max}]$ and it is clear that $\eta \le x \le \xi$ and $\eta\le y\le
  \xi$. Moreover, since all the norms are equivalent in $\real^n$,
  there exists $M_1,M_2>0$ such that $M_1\|v\| \le \|v\|_{\infty}\le
  M_2\|v\|$. As a result, we get
  $\|\eta-\xi\|\le M^{-1}_1\|\eta-\xi\|_{\infty}=M^{-1}_1
  \|x_0-y_0\|_{\infty} \le M^{-1}_1M_2\|x_0-y_0\|$. We set
  $M^{\frac{1}{2}} = M_1^{-1}M_2>0$. Since the vector field $f$ is monotone, for every $t\ge t_0$,
  \begin{align*}
    \phi(t,t_0,\eta)&\le \phi(t,t_0,x) \le \phi(t,t_0,\xi),\\
    \phi(t,t_0,\eta)&\le \phi(t,t_0,y) \le \phi(t,t_0,\xi).
  \end{align*}
  This means that
  \begin{align*}
    \|\phi(t,t_0,x)-\phi(t,t_0,y)\| & \le
    M^{-1}_1\|\phi(t,t_0,x)-\phi(t,t_0,y)\|_{\infty} \\ & \le M^{-1}_1
  \|\phi(t,t_0,\xi)-\phi(t,t_0,\eta)\|_{\infty}  \\ & \le  M^{-1}_1M_2
    \|\phi(t,t_0,\xi)-\phi(t,t_0,\eta)\| \\ & = M^{\frac{1}{2}} 
  \|\phi(t,t_0,\xi)-\phi(t,t_0,\eta)\|.
    \end{align*}
    However, we know that $\eta\le \xi$, thus by the above argument we get
    $\verti{\phi(t,t_0,\eta)-\phi(t,t_0,\xi)} \leq e^{b(t-s)}
    \verti{\phi(s,t_0,\eta)-\phi(s,t_0,\xi)}$. Therefore,
    \begin{align*}
      \|\phi(t,t_0,x_0)&-\phi(t,t_0,y_0)\|\le M^{\frac{1}{2}}
      \verti{\phi(t,t_0,\eta)-\phi(t,t_0,\xi)} \\ &
      \leq
      M^{\frac{1}{2}}  e^{b(t-s)} \|\phi(s,t_0,\eta)-\phi(s,t_0,\xi)\|
      \\ &= M
           e^{b(t-s)}\|\phi(s,t_0,x_0)-\phi(s,t_0,y_0)\|. \end{align*}\end{proof}

       \begin{remark} For monotonic norms, using
           Theorem~\ref{lem:lemma}\ref{p2:Metzler}, the notion of
           conic matrix measure coincides with the standard matrix
           measure on Metzler matrices. Therefore,
           Theorem~\ref{thm:contraction-monotone-general} can be
           completely recovered from~\cite[Theorem
           31]{AD-SJ-FB:20o}. However, for non-monotonic norms,
           Theorem~\ref{thm:contraction-monotone-general} provides a
           conic matrix measure condition for incremental exponential
           stability of the system. The next example elaborates this
           point in more detail.
         \end{remark}

       \begin{example}
           Consider the class of dynamical system on $\real^2$
           \begin{align}\label{eq:ex1}
             \dot{x}_1 &= -x_1 + \alpha x_2 -\gamma g(x_1),\nonumber\\
             \dot{x}_2 &= \beta x_1 - x_2,
           \end{align}
           where $\alpha,\beta,\gamma \ge 0$ and $g:\real\to
           \real_{\ge 0}$ is differentiable with $g(0)=0$ and $0\le g'(x) \le \overline{G}$, for every
           $x\in \real$. It is easy to see that, since
           $\alpha,\beta\ge 0$, the dynamical system~\eqref{eq:ex1} is
           monotone. For $\alpha =\frac{1}{2}$, $\beta=1.2$, $\gamma =
           1$, and $\overline{G}=0.1$, we have
           \begin{align*}
             Df(x_1,x_2) =\begin{bmatrix}
               -1 -g'(x_1) & \frac{1}{2}\\
               1.2 & -1
             \end{bmatrix}
                     \end{align*}
                     For monotonic norm
                     $\|\cdot\|_{\infty}$, we have $\mu^+_{\infty}(Df(x_1,x_2))
                     =\mu_{\infty}(Df(x_1,x_2)) = 0.2$. Therefore the
                     dynamical system~\eqref{eq:ex1} is not contracting with respect to $\ell_{\infty}$-norm. For the
                     non-monotonic norm $\|\cdot\|_{\infty,R}$ with
                     $R=\begin{bmatrix}-1 & 1 \\ 1 & 1\end{bmatrix}$,
                     we can compute
\begin{align}\label{eq:11}
                                                       \mu^+_{\infty,R}(Df(x_1,x_2))
                                                       & =
                                                       \mu^+_{\infty}(R
                                                       Df(x_1,x_2)
                                                       R^{-1}) \nonumber\\ & =
                                                       \mu^+_{\infty}\left(\begin{bmatrix}
                                                           -1.85
                                                           -\tfrac{g'(x_1)}{2}
                                                           & 0.35 +
                                                           \tfrac{g'(x_1)}{2}\\
                                                           -0.35+\tfrac{g'(x_1)}{2}
                                                           &
                                                           -0.15-\tfrac{g'(x_1)}{2}\end{bmatrix}\right)
                                                             \nonumber\\
                                                       & = -0.15 -
                                                         \tfrac{g'(x)}{2}
                                                     \le -0.15
\end{align}
Therefore, using Theorem~\ref{thm:contraction-monotone-general} and inequality~\eqref{eq:11}, the
dynamical system~\eqref{eq:ex1} is incrementally exponentially
stable with respect to the non-monotonic norm
$\|\cdot\|_{\infty,R}$. As a consequence, $\vect{0}_2$ is the globally
exponentially stable equilibrium point of the dynamical system~\eqref{eq:ex1}.  On the other hand, we have    
\begin{align}\label{eq:22}
  \mu_{\infty,R}(Df(x_1,x_2))
  & =
    \mu_{\infty}(R
    Df(x_1,x_2)
    R^{-1}) \nonumber\\ & =
                          \mu_{\infty}\left(\begin{bmatrix}
                              -1.85
                              -\tfrac{g'(x_1)}{2}
                              & 0.35 +
                              \tfrac{g'(x_1)}{2}\\
                              -0.35+\tfrac{g'(x_1)}{2}
                              &
                              -0.15-\tfrac{g'(x_1)}{2}\end{bmatrix}\right)
                                \nonumber\\ & \ge
                                              0.2
                                              -g'(x_1)
                                              \ge 0.1.
\end{align}
However, using~\cite[Theorem 31]{AD-SJ-FB:20o} and inequality~\eqref{eq:22}, the
dynamical system~\eqref{eq:ex1} is not contracting with respect to
$\|\cdot\|_{\infty,R}$.\end{example}

     We can also simplify Theorem~\ref{thm:contraction-monotone-general} for diagonally-weighted norms. 

\begin{corollary}[Diagonally-weighted norms]\label{thm:diagonal-l1-monotone}
  Let $\dot{x}=f(t,x)$ be a monotone dynamical
  system and $\eta\in \real^{n}_{>0}$. Then the following statements
  about $[\eta]$-weighted $\ell_1$-norm are equivalent:
  \begin{enumerate}
  \item\label{p1:1-norm} $\mu_{1,[\eta]}(\jac{f}(t,x)) \le b$, for every
    $(t,x)\in \real_{\ge 0}\times \real^n$;
  \item\label{p3:1-norm} $\eta^{\top}(f(t,x)-f(t,y))\le b\eta^{\top} (x-y)$, for every
    $x\ge y$ and every $t\in \real_{\ge 0}$;
    \item\label{p4:1-norm} $\|\phi(t,t_0,x_0)-\phi(t,t_0,y_0)\|_{1,[\eta]} \le e^{b(t-s)}
      \|\phi(s,t_0,x_0)-\phi(s,t_0,y_0)\|_{1,[\eta]}$, for every $x_0,y_0\in
      \real^n$ and every $t_0\le t\le s$. 
    \end{enumerate}
    Similarly, the following statements about $[\eta]^{-1}$-weighted $\ell_{\infty}$-norm are equivalent:
    \begin{enumerate}\setcounter{enumi}{3}
  \item\label{p1:inf-measure} $\mu_{\infty,[\eta]^{-1}}(\jac{f}(t,x)) \le b$, for every
    $(t,x)\in \real_{\ge 0}\times \real^n$;
  \item\label{p3:inf-integral} $f(t,x)-f(t,y)\le b(x-y)$, for every $t\in \real_{\ge 0}$ and every
    $x=y+ c \eta$ with $c>0$;
    \item\label{p4:inf-traj} $\|\phi(t,t_0,x_0)-\phi(t,t_0,y_0)\|_{\infty,[\eta]^{-1}} \le e^{b(t-s)}
      \|\phi(s,t_0x_0)-\phi(s,t_0,y_0)\|_{\infty,[\eta]^{-1}}$, for every $x_0,y_0\in
      \real^n$ and $t_0\le s\le t$. 
    \end{enumerate}
  \end{corollary}
  \smallskip
  \begin{proof} Note that $\|\cdot\|_{1,[\eta]}$ is a monotonic norm
    and its associated WP is given by~\eqref{eq:1-R} with
    $R=[\eta]$. Similarly, $\|\cdot\|_{\infty,[\eta]^{-1}}$ is a
    monotonic norm and its associated WP is given by~\eqref{eq:inf-R}
    with $R=[\eta]^{-1}$. Regarding
    \ref{p1:1-norm}~$\implies$\ref{p3:1-norm},
\begin{align*}
  \eta^{\top}&(f(t,x)-f(t,y))= \eta^{\top}\Big(\int_{0}^{1} Df(t,\tau x +
                               (1-\tau)y)d\tau\Big)(x-y)\\ & = \int_{0}^{1} \eta^{\top}Df(t,\tau x +
                                                               (1-\tau)y)(x-y) d\tau \le b\eta^{\top}(x-y),
\end{align*}
where the first equality is by the Mean Value Theorem and the last
inequality holds by Theorem~\ref{thm:computational}\ref{p1:l1-identity} and using the fact that $x-y\ge
\vect{0}_n$. Regarding~\ref{p3:1-norm}$\implies$\ref{p1:1-norm}, pick
$x = y + hv$, for $v\in \real^n_{\ge 0}$ and $h>0$. Thus, we get
\begin{align*}
  \eta^{\top}(f(t,x)-f(t,y))= \eta^{\top}(f(t,y+hv)-f(t,y))\le bh\eta^{\top}v.
\end{align*}
By taking the limit as $h\to 0^{+}$, for every $y\in \real^n$ and
every $v\in \real^n_{\ge 0}$,
\begin{align*}
  \eta^{\top}Df(t,y)v = \lim_{h\to 0^+}\eta^{\top}\tfrac{f(t,y+hv)-f(t,y)}{h}  \le
  b\eta^{\top}v. 
\end{align*}
The result then follows by
Theorem~\ref{thm:computational}\ref{p1:l1-identity}.

Regarding~\ref{p1:inf-measure} $\implies$\ref{p3:inf-integral}, note
that, for every $c>0$ such that $x = y + c\eta$, we have
\begin{align*}
  f(t,x) - f(t,y) &= \int_{0}^{1}
                    \jac{f}(t, (1-\tau)y + \tau x) (x-y) d\tau  \\ & = \int_{0}^{1}
                                                                     \jac{f}(t, y + \tau c\eta ) (c\eta) d\tau  \le b (x-y) 
\end{align*}
where the inequality follows from
Theorem~\ref{thm:computational}\ref{p2:linf-identity}. For~\ref{p3:inf-integral}
$\implies$\ref{p1:inf-measure},
\begin{align*}
  \jac{f}(t,x)\eta = \lim_{h\to 0^{+}}
  \frac{f(t,x+h\eta)-f(t,x)}{h} \le b\eta. 
\end{align*}
The result follows by
Theorem~\ref{thm:computational}\ref{p2:linf-identity}. The rest of the
proof follows from Theorem~\ref{thm:contraction-monotone-general}.
\end{proof}

   Next, we use the notion of conic matrix measure and weak pairing to study 
  exponential convergence of positive systems to their equilibrium points. 
   
   \begin{theorem}[Converging positive systems]\label{thm:positive}
     Let $\dot{x}=f(t,x)$ be a positive system with equilibrium point
     $\vect{0}_n$, $\|\cdot\|$ be a norm with associated WP
     $\WP{\cdot}{\cdot}$ satisfying Deimling's inequality and the curve
     norm derivative formula, and $b\in \real$. Consider 
      \begin{enumerate}[label=\textup{(A\arabic*)}]
      \item\label{p3:integral-general-0} $\WP{f(t,x)}{x}\le b\|x\|^2$,
        for every $(t,x)\in \real_{\ge 0}\times \real^n_{\ge 0}$;
     \item\label{p5:contraction-0}
        $\|\phi(t,t_0,x_0)\| \le
        e^{b(t-s)}\|\phi(s,t_0,x_0)\|$, for every $t_0\le s\le t$ and every
        $x_0\in \real^{n}_{\ge 0}$.
        \item \label{p6:diff-0} $\mu^+(B(t,x))\le b$, for
         every $x\ge \vect{0}_n$ and every $t\in \real_{\ge 0}$ where
         $B(t,x)$ is such that $f(t,x)=B(t,x)x$. 
     \end{enumerate}
\    Then the following statements hold: 
     \begin{enumerate}
     \item\label{p1:positive} \ref{p3:integral-general-0} 
       and~\ref{p5:contraction-0} are equivalent;
     \item\label{p2:positive} \ref{p6:diff-0} implies \ref{p3:integral-general-0} 
       and~\ref{p5:contraction-0}.
       \end{enumerate}
     \end{theorem}
     \smallskip
     \begin{proof}
       Regarding~\ref{p3:integral-general-0}$\iff$
       \ref{p5:contraction-0},
       the proof is similar to~\cite[Theorem 33]{AD-SJ-FB:20o} and we omit
       it.

       Regarding~\ref{p6:diff-0}$\implies$\ref{p5:contraction-0}, for every $x\ge 0$, we have
     \begin{align*}
     \WP{f(t,x)}{x} =
       \WP{B(t,x)x}{x} \le
       \mu^+(B(t,x))\|x\|^2\le b\|x\|^2.
    \end{align*}
    where the second inequality holds by Theorem~\ref{lem:lemma}\ref{p3:Lumer}. 
  \end{proof}

      In the next example, we investigate the role of conic matrix measures
      in the sufficient condition for exponential stability of
      positive systems in Theorem~\ref{thm:positive}. 

\begin{example}
  Consider the following dynamical system on $\real^2$:
  \begin{align}\label{ex:eq-contract}
    \dot{x}_1 &= -2x_1 + x_2:=f_1(x_1,x_2),\nonumber\\
    \dot{x}_2 &= - x_1\alpha(x_2) - x_2:=f_2(x_1,x_2),
  \end{align}
  where $\alpha:\real\to \real_{\ge 0}$ is a non-negative
  non-decreasing function.  First note that,
  $Df(x_1,x_2) = \left[\begin{smallmatrix} -2 & 1\\-\alpha(x_2) &
      -x_1\alpha'(x_2)-1\end{smallmatrix}\right]$.  Thus, the vector
  field $f$ is not monotone
  on $\real^2_{\ge 0}$ because $-\alpha(r)\le 0$ for every
  $r\in \real$. However, the dynamical system~\eqref{ex:eq-contract} is
  positive with an equilibrium point at
  $\vect{0}_2\in \real^2_{\ge 0}$. Moreover, 
  \begin{align*}
    f(x_1,x_2) = \begin{bmatrix}-2 & 1\\ -\alpha(x_2)
      & -1\end{bmatrix}\begin{bmatrix}x_1\\ x_2\end{bmatrix} :=
    B(x_1,x_2) \begin{bmatrix}x_1\\ x_2\end{bmatrix}.
  \end{align*}
  Using Theorem~\ref{thm:computational}\ref{p3:conic-inf}, for every $(x_1,x_2)\in \real^2_{\ge 0}$, 
  \begin{align}\label{eq:bounds}
    \mu_{\infty}^{+}(B(x_1,x_2))
  = -1 \le \mu_{\infty}(B(x_1,x_2)) = -1+\alpha(x_2).
  \end{align}
  By Theorem~\ref{thm:positive}\ref{p2:positive}, every trajectory $t\mapsto [x_1(t),x_2(t)]^{\top}$
  of the positive system~\eqref{ex:eq-contract} {\color{black}starting at
  $[x_1(0),x_2(0)]^{\top}\in \real^2_{\ge 0}$} satisfies
  \begin{align*}
    \|[x_1(t),x_2(t)]^{\top}\|_{\infty}\le e^{-t}
    \|[x_1(0),x_2(0)]^{\top}\|_{\infty}. 
  \end{align*}
  It is
  worth mentioning that, by equation~\eqref{eq:bounds}, the
  $\ell_{\infty}$-matrix measure of $B(x_1,x_2)$ might not be bounded
  and cannot be used to deduce 
  convergence of trajectories of~\eqref{ex:eq-contract} to $\vect{0}_2$.
  \end{example}

\section{Applications}

In this section, we present two applications for our
  non-Euclidean contraction framework for monotone and positive
  systems. As a first application, we show that a Hopfield neural
  network with excitatory interactions between its neurons is monotone
  but non-positive. We then use our framework to analyze stability and
  robustness of excitatory Hopfield neural networks. As a second
  application, we develop a framework for stability analysis of
  networks of interconnected systems using positive but non-monotone
  comparison systems.

\subsection{Excitatory Hopfield neural networks}

Hopfield model is a class of recurrent neural networks that can serve
as an associative memory system~\cite{JJH:1982}. There has been a
recent growing interest in the machine learning community to use
variations of Hopfield neural networks to store information or to
learn prototypes~\cite{HR-etal:20}. However, neural networks are
notoriously vulnerable to adversarial perturbations of their input;
small changes in their input can cause a large change in their
output~\cite{CZ-WZ-IS-JB-DE-IG-RF:13}. In this section, we study
stability and input-output robustness of Hopfield neural network with
excitatory neuron interactions and provide explicit adversarial
robustness guarantees for this class of learning algorithms. The
dynamics of the Hopfield neural network is given by
\begin{align}\label{eq:neural_network}
  \dot{x} = -\Lambda x + Tg(x) + I(t) : = F_{\mathrm{H}}(x),
\end{align}
where $x \in \real^n$ is the state of neurons, $\Lambda\in \real^{n\times
  n}$ is the diagonal positive-definite matrix of dissipation rates, $T \in
\real^{n\times n}$ is the interaction matrix, and $I:\real_{\ge 0}\to
\real^n_{\ge 0}$ is a time-varying input.
Assume $g(x) = (g_1(x_1),\ldots,g_n(x_n))^{\top}$, where the $i$th
activation function $g_i$ is Lipschitz continuous, monotonic non-decreasing
with $g_i(0)=0$ and with the finite sector property:
\begin{align*}
  0 \le \tfrac{g_i(x)-g_i(y)}{x-y}:=G_i(x,y) \le \overline{G}_i,
\end{align*}
where $\overline{G}= (\overline{G}_1,\ldots, \overline{G}_n)^{\top}\in
\real^{n}_{>0}$. We study excitatory Hopfield networks, i.e., neural
networks with Metzler interaction matrix $T$.

\begin{proposition}[Contracting Hopfield neural
  networks]\label{thm:neural-network} Consider the Hopfield neural
  network~\eqref{eq:neural_network} with an irreducible {\color{black}non-negative}
  interaction matrix
  $T$. Assume the Metzler matrix $-\Lambda+ T\overline{G}$ is Hurwitz with
  $(-c,v)$ and $(-c,w)$ its left and the right Perron eigenpair,
  respectively. For any $p\in[1,\infty]$, define $q\in [1,\infty]$ by
  $\frac{1}{p}+\frac{1}{q}=1$ (with convention $1/\infty = 0$) and $\eta\in
  \real^n_{>0}$ by
  \begin{align*}
    \eta =\Big(v_1^{\frac{1}{p}}/w_1^{\frac{1}{q}},\ldots,
    v_n^{\frac{1}{p}}/w_n^{\frac{1}{q}}\Big)^{\top}.
    \end{align*}
   Then the following statements hold for any $p\in[1,\infty]$:
  \begin{enumerate}
  \item\label{p2:contracting} the Hopfield neural
    network~\eqref{eq:neural_network} is monotone and contracting with respect
    to the norm $\|\cdot\|_{p,[\eta]}$ with rate $c$;
  \item\label{p3:constantinput} if $I(t) = I^*$ is constant, then the
    Hopfield neural network~\eqref{eq:neural_network} has a unique
    globally exponentially stable equilibrium point $x_I^*$ with the
    Lyapunov functions $\|x-x_I^*\|_{p,[\eta]}$ and
    $\|F_{\mathrm{H}}(x)\|_{p,[\eta]}$;

    \item\label{p5:input-output} if $t\mapsto x_I(t)$ and $t\mapsto
      x_J(t)$ are solutions of the Hopfield neural network~\eqref{eq:neural_network} for input signals $t\mapsto
      I(t)$ and $t\mapsto J(t)$ respectively, then, for every $t\in \real_{\ge 0}$,
      \begin{align*}
        \|x_I(t)-x_J(t)\|_{p,[\eta]} &\le
        e^{-ct}\|x_I(0)-x_J(0)\|_{p,[\eta]} \\ & + \int_{0}^{t} e^{-c(t-s)}\|I(s)-J(s)\|_{p,[\eta]}ds.
\end{align*}
  \end{enumerate}
\end{proposition}
\smallskip
\begin{proof}
 Let $i\in\{1,\ldots,n\}$ and consider $x\le y$
 such that $x_i=y_i$. For every $i\ne j$, by the finite sector
 property of $g_i$, we have $g_i(x_i)\le g_i(y_j)$, and thus
 \begin{align*}
   [F_{\mathrm{H}}(x)]_i & = -\gamma_ix_i + \sum\nolimits_{j=1}^{n} T_{ij}g_j(x_j) +
   I(t) \\ & \le -\gamma_iy_i + \sum\nolimits_{j=1}^{n} T_{ij}g(y_j) + I(t)= [F_{\mathrm{H}}(y)]_i,
 \end{align*}
 where the inequality holds because the matrix $T$ is Metzler. This means that the Hopfield
  neural network~\eqref{eq:neural_network} is monotone. Moreover, for every $x\ge
  y\ge \vect{0}_n$,
  \begin{align*}
    \|x-y\|_{p,[\eta]}&D^{+}\|x-y\|_{p,[\eta]} \\ & =\WP{-\Lambda(x-y) + T(g(x)-g(y))}{x-y}_{p,[\eta]}\\
    &\le \WP{(-\Lambda +  T\overline{G})(x-y)}{x-y}_{p,[\eta]} \\  
    &\le \mu_{p,[\eta]}(-\Lambda +    T\overline{G})\|x-y\|_{p,[\eta]}^2
    = -c \|x-y\|^2 _{p,[\eta]},
  \end{align*}
  where the first equality is the curve norm derivative formula. Since $g$
  is non-decreasing and $T$ is non-negative, we get the bound
  $T(g(x)-g(y))\le T\overline{G}(x-y)$ for every $x\ge y\ge
  \vect{0}_n$. Lemma~\ref{lem:useful}\ref{p2:inequalityWP} and this bound
  give us the second inequality. The third inequality holds by definition
  of $\mu_{p,[\eta]}$ and the fourth equality holds by~\cite{JA:96} using
  the fact that $-\Lambda + T\overline{G}$ is Metzler, irreducible, and
  Hurwitz with Perron eigenvalue $-c$ and left and right Perron
  eigenvectors $v$ and $w$. Then, parts~\ref{p2:contracting}
  and~\ref{p3:constantinput} follow from
  Theorem~\ref{thm:contraction-monotone-general}. Regarding
    part~\ref{p5:input-output}, 
    \begin{align*}
      \|x_{I}-&x_J\|_{p,[\eta]} D^+\|x_{I}-x_J\|_{p,[\eta]}
     \\ & = \WP{F_{\mathrm{H}}(x_I)-F_{\mathrm{H}}(x_J) +
      I(t)-J(t)}{x_I-x_J}_{p,[\eta]} \\ & \le
      \WP{F_{\mathrm{H}}(x_I)-F_{\mathrm{H}}(x_J)}{x_I-x_J}_{p,[\eta]}
      \\ & \qquad\qquad \qquad\qquad \quad+
      \WP{I(t)-J(t)}{x_I(t)-x_J(t)}_{p,[\eta]} \\ & \le
      -c\|x_I-x_J\|^2_{p,[\eta]}+
                                                    \|x_I-x_J\|_{p,[\eta]}\|I(t)-J(t)\|_{p,[\eta]},
    \end{align*}
    where the first equality holds by the curve norm derivative
    formula, the second inequality holds by subadditive property of
    WPs, and the third inequality holds by contractivity of
    $F_{\mathrm{H}}$ and the Cauchy-Schwarz inequality. 
    This implies that $D^+\|x_{I}(t)-x_J(t)\|_{p,[\eta]}\le
    -c\|x_I(t)-x_J(t)\|_{p,[\eta]}+\|I(t)-J(t)\|_{p,[\eta]}$, for
    every $t\in \real_{\ge 0}$. The
    result follows by Gr\"{o}nwall\textendash{}Bellman
    inequality~\cite[Lemma 11]{AD-SJ-FB:20o}.

\end{proof}

\begin{remark}[Comparison with the literature]
  We refer to~\cite{HZ-ZW-DL:14} for a review of stability properties of
  Hopfield neural networks; e.g., it is known that Hurwitzness of
  $-\Lambda+T\overline{G}$ (as we assume in
  Proposition~\ref{thm:neural-network}) implies global exponential
  stability.  To the best of our knowledge, the strong contractivity (with
  respect to appropriately weighted $p$-norms) in part~\ref{p2:contracting}
  and the Lyapunov functions in part~\ref{p3:constantinput} are novel.
  Recall that, as reviewed in the Introduction, strong contractivity is a
  stronger property than global exponential stability. To the best of our knowledge, the
  input-to-output stability in part~\ref{p5:input-output} is novel and is directly
  applicable to obtain adversarial robustness guarantees of Hopfield neural
  networks. \oprocend
  \end{remark}

  \subsection{Non-monotone comparison systems}

    Comparison principles are well-established techniques in
    dynamical system theory to infer stability of a
    dynamical system using properties of a simpler comparison
    systems. In most of the existing
    comparison frameworks in the literature, monotonicity of the
    comparison system plays a crucial role~\cite{AAM-AYO:90,SGN-WMH:06}. In this subsection, we
    develop a novel comparison principle for stability analysis of
    networks of interconnected systems. Unlike the existing comparison
    results in the literature, our framework uses positive
    comparison systems which are not necessarily monotone. Consider
  the interconnection of $n$ subsystems:
  \begin{align}\label{eq:interconnection}
    \dot{x}_i= f_i(x,u_i), \qquad i\in \{1,\ldots,n\}.
  \end{align}
  where $x_i\in \real^{N_i}$ is the state and $u_i\in \real^{M_i}$ is the
  exogenous input for the $i$th subsystem. We define $N=\sum_{i=1}^{n} N_i$
  and $M =\sum_{i=1}^{n} M_i$ and $x=(x_1,\ldots,x_n)^{\top}\in
  \real^{N}$. We assume that, for every $i\in \{1,\ldots,n\}$, we have
  $f_i(\vect{0}_N,\vect{0}_{M_i}) = \vect{0}_{N_i}$ and the $i$th subsystem
  is equipped with the norm $\|\cdot\|_i$ on $\real^{N_i}$. We assume
  that, for every $i\in \{1,\ldots,n\}$, the $i$th subsystem has a storage function $V_{i}:\real^{N_i}\to \real_{\ge 0}$ such that
 \begin{enumerate}[label=\textup{(B\arabic*)}]
  \item\label{con1:lyap} there exists class $\mathcal{K}_{\infty}$
    functions $\underline{\alpha}_i$ and $\overline{\alpha}_i$ such
    that $\underline{\alpha}_i(\|x_i\|_i) \le V_i(x_i) \le
    \overline{\alpha}_i(\|x_i\|_i)$;
    \item\label{con3:compareV} for every $x\in \real^N,u\in \real^M$, the inequality
         \begin{align}\label{eq:dissipative-inequality}
           \mathcal{L}_{f_i} V_i(x_i) \le -\alpha_i(V_i(x_i))  + g_i(V(x)) +
           \gamma_i(u).
         \end{align}
         holds for some function $g_i:\real^{N}_{\ge 0}\to
         \real$ with $g_i(\vect{0}_{N})=0$ and for
         some class $\mathcal{K}_{\infty}$ functions $\alpha_i$ and
         $\gamma_i$.
       \end{enumerate}
       We define the maps $\Gamma:\real^N_{\ge 0}\to \real^n$ and
       $A:\real^N_{\ge 0}\to \real^n$ by:
       \begin{align*}
         \Gamma(x)=(g_1(x),\ldots,g_n(x))^{\top},\quad A(x)=(\alpha_1(x_1),\ldots,\alpha_n(x_n))^{\top}
       \end{align*}
       One can also define the \emph{non-monotone comparison system} by
       \begin{align}\label{eq:comparisonsystem}
         \dot{v}_i = -A_i(v_i) + \Gamma_i(v) +
         \gamma_{i}(u),\;\;\;  i\in \{1,\ldots,n\}. 
       \end{align}
       Using the inequality~\eqref{eq:dissipative-inequality}, one can
       show that, for $u=\vect{0}_M$, the comparison system~\eqref{eq:comparisonsystem} is
       a positive dynamical system. However, since $g_i$ can be any
       arbitrary function, the comparison system~\eqref{eq:comparisonsystem} is not
       necessarily monotone.  
       \begin{remark}[Input-to-state stability]    
         The inequality \eqref{eq:dissipative-inequality} can be
         considered as a generalization of component-wise
         input-to-state stability (ISS). An interconnected
         system is component-wise ISS if each of its subsystems is ISS
         when interconnections between subsystems are considered as
         the input. In other words, the interconnected
         system~\eqref{eq:interconnection} is component-wise ISS if,
         for every $i\in \{1,\ldots,n\}$, the storage function $V_i$ satisfies
         \begin{align*}
           \mathcal{L}_{f_i} V_i(x_i) \le -\alpha_i(V_i(x_i)) + \sum\nolimits_{j\ne i} \gamma_{ij}(V_j(x_j)) + \gamma_{iu}(u_i)
         \end{align*}
         for class $\mathcal{K}_{\infty}$ function $\alpha_i$ and
         class $\mathcal{K}$ functions $\gamma_{ij}$ and
         $\gamma_{iu}$. Indeed, if the interconnected
         system is component-wise ISS, then the associated comparison system~\eqref{eq:comparisonsystem}
         is monotone. 
       \end{remark}

 \begin{proposition}[Stability of interconnection of systems]\label{thm:ISS-interconnection-1}
   Consider the interconnected
   system~\eqref{eq:interconnection} and suppose that every subsystem
   satisfies conditions~\ref{con1:lyap} and~\ref{con3:compareV}
   above. Let $p\in [1,\infty]$ and $R\in \real^{n\times n}$ be an
   invertible non-negative matrix. Suppose that there exists $c>0$
   such that, for every $v\ge
     \vect{0}_n$, 
   \begin{align}
  \label{cond3:integral-dem}
       -\WP{-A(v)}{v}_{p,R}\ge
     \WP{\Gamma(v)}{v}_{p,R}+c\|v\|_{p,R}^2 \tag{$C_1$}
     \end{align}
     Then, the following statements hold:
  \begin{enumerate}
  \item\label{p1.1:comparison_system} the comparison
    system~\eqref{eq:comparisonsystem} converges
    exponentially to $\vect{0}_n$
  \item\label{p1:0-converge} for $u(t)=\vect{0}_M$, every trajectory
    of the interconnected system~\eqref{eq:interconnection} converges
    to $\vect{0}_N$.
 
     \item\label{p2:ISS-zero} the system~\eqref{eq:interconnection} is
    input-to-state stable in the sense that, for every $i\in
    \{1,\ldots,n\}$ and $t\ge 0$, there exists $L_i>0$, such that
    \begin{multline}\label{p2:ISS}
       \|x_i(t)\|_i \le
       \underline{\alpha}_i^{-1}\big(L_i e^{-{c}t}\|V(x(0))\|_{p,R}+\\ 
       \tfrac{L_i(1-e^{-{c}t})}{c}\max_{\tau\in [0,t]}\|\gamma(u(\tau))\|_{p,R}\big).
     \end{multline}
    \end{enumerate}
     Alternatively, if $v\mapsto -A(v) + \Gamma(v)$ is continuously differentiable,
     then~\ref{p1.1:comparison_system}, \ref{p1:0-converge}, and~\ref{p2:ISS-zero} still holds by replacing
     condition~\eqref{cond3:integral-dem} with the following stronger
     condition:
\begin{align}\label{cond2:monotone}
         \mu _{p,R}^{+}(B(v))\le -c \tag{$C_2$}
\end{align}
where $B(v)\in \real^{n\times n}$ satisfies $B(v)v = -A(v) + \Gamma(v)$, for $v\in \real^n_{\ge 0}$.
\end{proposition}
     \smallskip
     \begin{proof}
       Regarding part~\ref{p1.1:comparison_system}, for $u=\vect{0}_M$, we have
       \begin{align*}
         \WP{-A(v) + \Gamma(v)}{v}_{p,R} \le \WP{-A(v)}{v}_{p,R} +\WP{\Gamma(v)}{v}_{p,R}  \le -c\|v\|^2_{p,R}
       \end{align*}
       for every $v\in \real^n_{\ge 0}$. Since the comparison system
       is positive, the result follows from Theorem~\ref{thm:positive}\ref{p1:positive}.  Regarding
    part~\ref{p1:0-converge}, by setting $V(x(t))=V(t)$, we get  
     \begin{align*}
       \|V(t)\|& D^{+}\|V(t)\|_{p,R} = \WP{\dot{V}(t)}{V(t)}_{p,R}\\
       &\le       \WP{-A(V(t)) + \Gamma(V(t)) + \gamma(u)}{V(t)}_{p,R}
       \\ & \le \WP{-A(V(t))}{V(t))}_{p,R} + 
       \WP{\Gamma(V(t))}{V(t)}_{p,R} \\ & \qquad\qquad \qquad\qquad \qquad\quad+  \WP{\gamma(u)}{V(t)}_{p,R}
       \\
       & \le -c
       \|V(t)\|_{p,R}^2 +  \|\gamma(u)\|\|V(t)\|_{p,R}.
     \end{align*}
     where the first equality holds by the curve norm derivative
     formula, the second inequality holds by Lemma~\ref{lem:useful}\ref{p2:inequalityWP}, the
     third inequality holds by subadditivity of the WPs, and
     the fourth inequality holds by the Cauchy-Schwarz
     inequality. This implies that
       \begin{align*}
         \|V(t)\|_{p,R} \le e^{-ct}\|V(0)\|_{p,R}+
       \tfrac{1-e^{-ct}}{c}\max_{\tau\in [0,t]}\|\gamma(u(\tau))\|_{p,R}.
       \end{align*}
       Therefore, for $u=\vect{0}_M$, we have
       $t\mapsto V(t)$ converges exponentially to $\vect{0}_n$ and thus
       $\lim_{t\to \infty} x(t) = \vect{0}_N$. Regarding part~\ref{p2:ISS-zero}, since $R$ is non-negative and invertible, there exists $L_i>0$
       such that
       $V_i(x_i) \le L_i\|V(x)\|_{p,R}$, for every $i\in
       \{1,\ldots,n\}$. Moreover, we know that
       $\underline{\alpha}_i(\|x_i\|_i) \le V_i(x_i)$, for every
       $i\in \{1,\ldots,n\}$. The result then easily
       follows. Finally, for continuously differentiable
       $v\mapsto -A(v)+\Gamma(v)$, condition~\eqref{cond2:monotone}
       implies condition~\eqref{cond3:integral-dem} by
       Theorem~\ref{thm:positive}\ref{p2:positive}.\end{proof}

     \begin{remark}[Small-gain interpretation]
       \begin{enumerate}     
  \item For condition~\eqref{cond3:integral-dem} the term
  $-\WP{-A(v)}{v}_{p,R}$ captures the incremental dissipation gains of the
  subsystems while the term $\WP{\Gamma(v)}{v}_{p,R}$ captures the
  incremental interconnection gains between subsystems. Therefore, one can
  interpret the condition \eqref{cond3:integral-dem} as a small-gain condition
  requiring the dissipation gains to dominate the interconnection gains.
      
\item For monotone vector field $\Gamma$, one can choose $p=1$ and
  $R=[\eta]\in \real^{n\times n}$ for some $\eta\in \real^n_{>0}$ and
  using Corollary~\ref{thm:diagonal-l1-monotone} to write condition
  \eqref{cond3:integral-dem} as:
  \begin{align*}
    \eta^{\top}A(v)\ge
    \eta^{\top}\Gamma(v)+ c\eta^{\top}v,
  \end{align*}
  for every $v\ge \vect{0}_n$. This result is similar to the
  small-gain theorem developed
  in~\cite{SND-BSR-FRW:10,BSR-CMK-SRW:10}.  
  
\item Compared to the classical comparison results
    (see~\cite{AAM-AYO:90,SGN-WMH:06,BSR-CMK-SRW:10}),
  Proposition~\ref{thm:ISS-interconnection-1} does not require
  monotonicity of the comparison system. Instead
  Proposition~\ref{thm:ISS-interconnection-1} is based on comparing the
  interconnected system with a positive comparison system. As a
  result, contrary to the existing small-gain theorems
  (see~\cite{SND-BSR-FRW:10,BSR-CMK-SRW:10}),
  Proposition~\ref{thm:ISS-interconnection-1} can take into account both
  the inhibitory and excitatory nature of the interactions between the
  subsystems. The next example illustrates this point in more
  detail.\oprocend
         \end{enumerate}
     \end{remark}

\begin{example}
         Consider the following system on $\real^2$:
  \begin{align}\label{ex:2}
    \dot{x}_1 &= -x_1 +\beta(x_2)x_1x_2^5 - 2x^3_1x_2^4\nonumber\\
    \dot{x}_2 &= -x_2 + x_1^6x_2-x_1^4x^3_2,
  \end{align}
  where $\beta:\real\to \real$ is such that {\color{black}$|\beta(r)| \le |r|$, for
  every $r\in \real$}. We choose the storage functions $V_i(x_i)=
  x_i^2$ for $i\in \{1,2\}$. One can construct a monotone comparison
  system for the dynamics~\eqref{ex:2} as follows:
  \begin{align*}
    \dot{V}_1 &=-2V_1+2\beta(x_2)x_2^5V_1 -4V^2_1V^2_2\le -2V_1+2V^3_2V_1, \\
    \dot{V}_2 &=-2V_2+2V^3_1V_2 -2V^2_1V_2^2\le -2V_2+ 2V^3_1V_2.
  \end{align*}
  Therefore, the comparison system has the form $\dot{v} = h(v)$ with
  $h(v) =\begin{bmatrix}-2v_1+2v^3_2v_1\\ -2v_2+
    2v^3_1v_2\end{bmatrix}$. Since the Jacobian of $h$
  is Metzler on $\real^2_{\ge 0}$, the comparison
  system $h$ is monotone on $\real^2_{\ge 0}$. However, this comparison system has two equilibrium points
  $v_1=v_2=0$ and $v_1=v_2=1$. Therefore, it is not possible to use
  comparison system $h$ to deduce
  global stability of $\vect{0}_2$ for the original dynamical system~\eqref{ex:2}.  On the
  other hand, one can construct a positive non-monotone comparison
  system for the dynamics~\eqref{ex:2} as follows:
  \begin{align*}
    \dot{V}_1 &=-2V_1+2\beta(x_2)x_2^5V_1 -4 V_1^2V^2_2\le
                -2V_1+2V_1V^3_2 -4 V_1^2V^2_2, \\
    \dot{V}_2 &=-2V_2+2V^3_1V_2 -2 V^2_1V_2^2\le -2V_2+2 V_2V^3_1 -2 V^2_1V_2^2.
  \end{align*}
  Therefore, the comparison system has the from $\dot{v}=A(v) +\Gamma(v)$ with
  $A(v_1,v_2)=\begin{bmatrix}-2v_1\\-2v_2\end{bmatrix}$ and
  $\Gamma(v_1,v_2) = \begin{bmatrix}2v_1v_2^3-4v_1^2v_2^2\\2 v_2v^3_1
    -2v^2_1v_2^2\end{bmatrix}$. We can also define $
    B(v_1,v_2)=2\begin{bmatrix}-1 &
    v_1v_2^2-2v_2v_1^2\\ v_2v_1^2-v_1v_2^2 & -1\end{bmatrix}$ where
  $B(v)v =-A(v)+\Gamma(v)$. Moreover, we get
  \begin{align*}
    \mu^{+}_{2}(B(v_1,v_2))
 & = \sup_{x\ge \vect{0}_2} \frac{x^{\top}Bx}{\|x\|^2_2} = 2\sup_{x\ge
    \vect{0}_2} \frac{x^{\top}\left[\begin{smallmatrix}-1 &-v_1^2v_2\\ 0 &
      -1\end{smallmatrix}\right]x}{\|x\|^2_2} \\ & =2\sup_{x\ge \vect{0}_2}
                                        \frac{-x_1^2-x_2^2-v_1^2v_2x_1x_2}{\|x\|^2_2}
                                        \le -2,
  \end{align*}
  where the first equality holds by Theorem~\ref{lem:lemma}\ref{p3:Lumer} and the second
  equality holds by the fact that $x^{\top}Bx =
  \tfrac{1}{2}x^{\top}(B+B^{\top})x$. Therefore,
  condition~\eqref{cond2:monotone} holds and, by Proposition~\ref{thm:ISS-interconnection-1}, every trajectory of the system~\eqref{ex:2} converges
to $\vect{0}_n$. 
       \end{example}

     \section{Conclusion}
     In this paper, we used conic matrix measures and weak pairings to
     characterize contracting monotone systems and to provided
     sufficient conditions for exponential convergence of positive
     systems to their equilibriums. As applications, we used our
     monotone contraction results to study contractivity and
     robustness of Hopfield neural networks. We also used our positive
     contraction results to established a novel and less-conservative
     framework for studying stability of interconnected networks. Future work
     includes extension of this framework to study monotone and
     positive systems which are weak- or
     semi-contracting~\cite{SJ-PCV-FB:19q} and to characterize
     contractivity of systems that are monotone with respect to
     arbitrary cones.

     \bibliographystyle{plain} \bibliography{alias,Main,FB,Ref,New}

\begin{thebibliography}{10}

\bibitem{JA:96}
J.~Albrecht.
\newblock Minimal norms of nonnegative irreducible matrices.
\newblock {\em Linear Algebra and its Applications}, 249(1):255--258, 1996.

\bibitem{ZA-EDS:14b}
Z.~Aminzare and E.~D. Sontag.
\newblock Contraction methods for nonlinear systems: {A} brief introduction and
  some open problems.
\newblock In {\em {IEEE} Conf.\ on Decision and Control}, pages 3835--3847,
  December 2014.

\bibitem{FB:22-CTDS}
F.~Bullo.
\newblock {\em Contraction Theory for Dynamical Systems}.
\newblock Kindle Direct Publishing, {1.0} edition, 2022.

\bibitem{GC-EL-KS:15}
G.~Como, E.~Lovisari, and K.~Savla.
\newblock Throughput optimality and overload behavior of dynamical flow
  networks under monotone distributed routing.
\newblock {\em IEEE Transactions on Control of Network Systems}, 2(1):57--67,
  2015.

\bibitem{SC:19}
S.~Coogan.
\newblock A contractive approach to separable {Lyapunov} functions for monotone
  systems.
\newblock {\em Automatica}, 106:349--357, 2019.

\bibitem{SND-BSR-FRW:10}
S.~N. Dashkovskiy, B.~S. R\"uffer, and F.~R. Wirth.
\newblock Small gain theorems for large scale systems and construction of {ISS}
  {Lyapunov} functions.
\newblock {\em SIAM Journal on Control and Optimization}, 48(6):4089--4118,
  2010.

\bibitem{AD-SJ-FB:20o}
A.~Davydov, S.~Jafarpour, and F.~Bullo.
\newblock {Non-Euclidean} contraction theory for robust nonlinear stability.
\newblock {\em IEEE Transactions on Automatic Control}, 2022.

\bibitem{MdB-DF-GR-FS:16}
M.~{Di~Bernardo}, D.~Fiore, G.~Russo, and F.~Scafuti.
\newblock Convergence, consensus and synchronization of complex networks via
  contraction theory.
\newblock In J.~L{\"u}, X.~Yu, G.~Chen, and W.~Yu, editors, {\em Complex
  Systems and Networks}, pages 313--339. Springer, 2016.

\bibitem{LF-SR:00}
L.~Farina and S.~Rinaldi.
\newblock {\em Positive Linear Systems: Theory and Applications}.
\newblock John Wiley \& Sons, 2000.

\bibitem{HRF-BB-MJ:18}
H.~R. {Feyzmahdavian}, B.~{Besselink}, and M.~{Johansson}.
\newblock Stability analysis of monotone systems via max-separable {Lyapunov}
  functions.
\newblock {\em IEEE Transactions on Automatic Control}, 63(3):643--656, 2018.

\bibitem{MG-SB:11-cvx}
M.~Grant and S.~Boyd.
\newblock {CVX}: Matlab software for disciplined convex programming, version
  2.1, March 2014.

\bibitem{MWH-HLS:03}
M.~W. Hirsch and H.~L. Smith.
\newblock Competitive and cooperative systems: A mini-review.
\newblock In L.~Benvenuti, A.~De Santis, and L.~Farina, editors, {\em Positive
  Systems}, pages 183--190, 2003.

\bibitem{JH-KS:98}
J.~Hofbauer and K.~Sigmund.
\newblock {\em Evolutionary Games and Population Dynamics}.
\newblock Cambridge University Press, 1998.

\bibitem{JJH:1982}
J.~J. Hopfield.
\newblock Neural networks and physical systems with emergent collective
  computational abilities.
\newblock {\em Proceedings of the National Academy of Sciences},
  79(8):2554--2558, 1982.

\bibitem{SJ-PCV-FB:19q}
S.~Jafarpour, P.~Cisneros-Velarde, and F.~Bullo.
\newblock Weak and semi-contraction for network systems and diffusively-coupled
  oscillators.
\newblock {\em IEEE Transactions on Automatic Control}, 67(3):1285--1300, 2022.

\bibitem{YK-BB-MC:20}
Y.~Kawano, B.~Besselink, and M.~Cao.
\newblock Contraction analysis of monotone systems via separable functions.
\newblock {\em IEEE Transactions on Automatic Control}, 65(8):3486--3501, 2020.

\bibitem{AK-TB-BG:16}
A.~Khanafer, T.~Ba{\c s}ar, and B.~Gharesifard.
\newblock Stability of epidemic models over directed graphs: {A} positive
  systems approach.
\newblock {\em Automatica}, 74:126--134, 2016.

\bibitem{WL-JJES:98}
W.~Lohmiller and J.-J.~E. Slotine.
\newblock On contraction analysis for non-linear systems.
\newblock {\em Automatica}, 34(6):683--696, 1998.

\bibitem{IRM-JJES:17}
I.~R. {Manchester} and J.-J.~E. Slotine.
\newblock Control contraction metrics: Convex and intrinsic criteria for
  nonlinear feedback design.
\newblock {\em IEEE Transactions on Automatic Control}, 62(6):3046--3053, 2017.

\bibitem{MM-EDS-TT:14}
M.~Margaliot, E.~D. Sontag, and T.~Tuller.
\newblock Entrainment to periodic initiation and transition rates in a
  computational model for gene translation.
\newblock {\em PLoS One}, 9(5):e96039, 2014.

\bibitem{AAM-AYO:90}
A.~A. Martynyuk and A.~Y. Obolensky.
\newblock On the theory of one-sided models in spaces with arbitrary cones.
\newblock {\em Journal of Applied Mathematics and Stochastic Analysis},
  3(2):85--97, 1990.

\bibitem{AGM:11}
A.~G. Mazko.
\newblock Positivity, robust stability and comparison of dynamic systems.
\newblock {\em Discrete and Continuous Dynamical Systems. Series A}, pages
  1042--1051, 2011.

\bibitem{SGN-WMH:06}
S.~G. Nersesov and W.~M. Haddad.
\newblock On the stability and control of nonlinear dynamical systems via
  vector lyapunov functions.
\newblock {\em IEEE Transactions on Automatic Control}, 51(2):203--215, 2006.

\bibitem{HR-etal:20}
H.~Ramsauer, B.~Schafl, J.~Lehner, P.~Seidl, M.~Widrich, T.~Adler, L.~Gruber,
  M.~Holzleitner, M.~Pavlovic, G.~K. Sandve, et~al.
\newblock Hopfield networks is all you need.
\newblock {\em arXiv preprint arXiv:2008.02217}, 2020.

\bibitem{AR:15}
A.~Rantzer.
\newblock Scalable control of positive systems.
\newblock {\em European Journal of Control}, 24:72--80, 2015.

\bibitem{BSR-CMK-SRW:10}
B.~S. R{\"u}ffer, C.~M. Kellett, and S.~R. Weller.
\newblock Connection between cooperative positive systems and integral
  input-to-state stability of large-scale systems.
\newblock {\em Automatica}, 46(6):1019--1027, 2010.

\bibitem{HLS:95}
H.~L. Smith.
\newblock {\em Monotone Dynamical Systems: An Introduction to the Theory of
  Competitive and Cooperative Systems}.
\newblock American Mathematical Society, 1995.

\bibitem{EDS:07}
E.~D. Sontag.
\newblock Monotone and near-monotone biochemical networks.
\newblock {\em Systems and Synthetic Biology}, 1(2):59--87, 2007.

\bibitem{CZ-WZ-IS-JB-DE-IG-RF:13}
C.~Szegedy, W.~Zaremba, I.~Sutskever, J.~Bruna, D.~Erhan, I.~Goodfellow, and
  R.~Fergus.
\newblock Intriguing properties of neural networks.
\newblock In {\em International Conference on Learning Representations}, 2014.

\bibitem{LV-SB:96}
L.~Vandenberghe and S.~Boyd.
\newblock Semidefinite programming.
\newblock {\em SIAM Review}, 38(1):49--95, March 1996.

\bibitem{HZ-ZW-DL:14}
H.~{Zhang}, Z.~{Wang}, and D.~{Liu}.
\newblock A comprehensive review of stability analysis of continuous-time
  recurrent neural networks.
\newblock {\em IEEE Transactions on Neural Networks and Learning Systems},
  25(7):1229--1262, 2014.

\end{thebibliography}
         
\end{document}